\documentclass{amsart}
\usepackage[T1]{fontenc}

\newtheorem{theorem}{Theorem}
\newtheorem{lemma}{Lemma}
\newtheorem{corollary}{Corollary}
\newtheorem{proposition}{Proposition}

\theoremstyle{definition}
\newtheorem{definition}{Definition}
\newtheorem{example}{Example}
\newtheorem*{notation}{Notation}

\theoremstyle{remark}

\theoremstyle{remark}
\newtheorem{claim}{Claim}

\newtheorem{remark}{Remark}

\newcommand{\baton}[1]{{\mathbb #1}}
\newcommand{\N}{\baton N}

\newcommand{\Z}{{\baton Z}}

\newcommand{\C}{{\baton C}}

\newcommand{\CF}{{\mathcal F}}
\newcommand{\CG}{{\mathcal G}}

\newcommand{\CP}{{\mathcal P}}
\newcommand{\CQ}{{\mathcal Q}}

\newcommand{\simp}{\underset{\CP}{\sim}}
\newcommand{\simqs}{\underset{\CP_s}{\sim}}

\newcommand{\equivp}{\underset{\CP}{\equiv}}
\newcommand{\equivq}{\underset{\CQ}{\equiv}}
\newcommand{\equivqy}{\underset{\CQ_Y}{\equiv}}

\newcommand{\approxq}{\underset{\CQ}{\approx}}

\newcommand{\one}{\boldsymbol{1}}

\newcommand{\inv}{^{-1}}
\DeclareMathOperator{\id}{\bf Id}
\newcommand{\wec}[1]{\langle#1\rangle}

\newcommand{\ov}[1]{\langle #1\rangle}
\newcommand{\br}[1]{[#1]}

\newcommand{\bx}{\mathbf{x}}
\newcommand{\by}{\mathbf{y}}
\newcommand{\bg}{\mathbf{g}}

\newcommand{\bh}{\mathbf{h}}
\newcommand{\bz}{\mathbf{z}}

\newcommand{\bu}{\mathbf{u}}
\newcommand{\bv}{\mathbf{v}}

\newcommand{\bgamma}{\boldsymbol{\gamma}}
\newcommand{\ubx}{\mathbf{\underline x}}
\newcommand{\uby}{\mathbf{\underline y}}

\newcommand{\ubu}{\mathbf{\underline u}}

\newcommand{\ubg}{\mathbf{\underline g}}
\newcommand{\ubh}{\mathbf{\underline h}}
\newcommand{\ubj}{\mathbf{\underline j}}

\newcommand{\type}[1]{^{[#1]}}
\newcommand{\norm}[1]{\lVert #1\rVert}

\DeclareMathOperator{\trace}{Trace}
\begin{document}

\title{Parallelepipeds, nilpotent groups, and Gowers norms}

\author{Bernard Host \and Bryna Kra}
\address{\'Equipe d'analyse et de math\'ematiques appliqu\'ees,
Universit\'e de Marne la Vall\'ee, 77454 Marne la Vall\'ee
Cedex, France}
\email{bernard.host@univ-mlv.fr}
\address{Department of Mathematics, Northwestern University,
2033 Sheridan Road,  Evanston, IL 60208-2730, USA}
\email{kra@math.northwestern.edu}
\thanks{The second author  was partially supported by 
NSF grant DMS-0555250.}
\date{\today}

\begin{abstract}
In his proof of 
Szemer\'edi's Theorem, Gowers introduced certain norms 
that are defined on a parallelepiped structure.  
A natural question is on which sets a parallelepiped
structure (and thus a Gowers norm) can be defined.  
We focus on dimensions $2$ and $3$ and show when 
this possible, and describe a correspondence
between the parallelepiped structures nilpotent groups.
\end{abstract}

\maketitle

\section{Introduction}
In his proof of Szemer\'edi's Theorem~\cite{S}, 
Gowers~\cite{gowers} introduced certain norms for functions on 
$\Z/N\Z$.  Shortly thereafter, these norms have been 
adapted for a variety of other uses, including Green 
and Tao's  major breakthrough showing that the primes contain 
arbitrarily long arithmetic progression~\cite{GT}.  
Similar objects were independently introduced by 
the authors and used to show convergence of some multiple 
ergodic averages~\cite{HK}, and since then have been used for a 
variety of related problems in ergodic theory.  
Our goal here is to introduce and describe the most 
general context in which the first two Gowers norms 
can be defined.  We call a `parallelogram structure,' respectively a
`parallelepiped structure,' the weakest structure a set 
must have so that one can define 
a Gowers $2$-norm, respectively a Gowers $3$-norm, 
on the set.

The first Gowers norm is the absolute value of the sum of 
the values of the function, and in fact is only a seminorm.  
The second Gowers norm can be completely described using Fourier 
analysis (in terms of the $\ell^4$ norm 
of the Fourier transform), and thus is closely linked to an abelian group 
on the circle.  Analogously, in ergodic theory 
the second seminorm can be characterized completely by the 
Kronecker factor in a measure preserving system (see 
Furstenberg~\cite{F1}), which is the largest 
abelian group rotation factor.  
The third Gowers norm is less well understood 
and can not  be simply described in terms of Fourier analysis.  
In ergodic theory the third seminorm corresponds to a $2$-step 
nilsystem, and more generally the $k$-th seminorm corresponds 
to a $(k-1)$-step nilsystem (see~\cite{HK}).  
In combinatorics, Green and Tao~\cite{GT2} have 
recently given a weak inverse theorem, but the 
correspondence with a $2$-step nilpotent group is not yet 
completely understood.  
We give conditions on a set that explain to what 
extent the correspondence with nilsystems can be made 
precise.

We start by defining the Gowers norms for $k\geq 2$.  
Let $\CP$ denote the subset
$$
 \{ (x_{00}, x_{01}, x_{10}, x_{11})\in (\Z/N\Z)^4\colon 
x_{00}-x_{01}-x_{10}+x_{11}=0\}
$$
of $(\Z/N\Z)^4$. 
For a function $f\colon \Z/N\Z\to\C$, 
the \emph{second Gowers norm $\|f\|_{U_2}$} is 
defined\footnote{The notation $\|\cdot\|_{U_k}$ was introduced later in 
the work of Green and Tao~\cite{GT}.}  
by
$$
\norm f_{U_2}^4 = 
\sum_{(x_{00}, x_{01}, x_{10}, x_{11})\in\CP}
f(x_{00})\overline{f(x_{01})}\, \overline{f(x_{10})}f(x_{11})\ .
$$
(Although this agrees with Gowers's original definition, Green 
and Tao prefer to normalize the sum and define the norm 
as the average of the sum.  In our context we prefer to 
work with the sum.)  
To define the norms  $U_k$  for $k\geq 3$, we need some notation.
\begin{notation}
The elements of $\{0,1\}^k$ are written without commas and 
parentheses. For $\epsilon=\epsilon_1\dots\epsilon_k\in\{0,1\}^k$ we 
write 
$$\lvert\epsilon\rvert = \epsilon_1+\dots+\epsilon_k\ .$$
For $x=(x_1,\dots,x_k)\in(\Z/N\Z)^k$ and $\epsilon\in\{0,1\}^k$, we 
write 
$$
 \epsilon\cdot x=\epsilon_1x_1+\dots+\epsilon_kx_k\ .
$$
$C\colon\C\to\C$ for complex conjugation. Therefore, for 
$n\in\N\cup\{0\}$ and $\xi\in\C$,
$$
 C^n\xi =
\begin{cases} \xi &\text{ if $n$ is even}\\
\overline\xi &\text{ if $n$ is odd\ .}
\end{cases}
$$
\end{notation}
For $k\geq 3$, the $U_k$ norm is the sum over $k$-dimensional 
parallelepipeds:
$$
\|f\|_{U_k}^{2^k} = \sum_{x\in (\Z/N\Z)^k}\sum_n
\prod_{\epsilon\in\{0,1\}^k} C^{|\epsilon|}f(n-\epsilon\cdot x)\ .
$$

For $k=3$, $\|f\|_{U_3}^{8}$
can be written as the sum 
\begin{multline}
\sum_{x,m,n,p\in \Z/N\Z} f(x)\overline{f(x+m)}\,
\overline{f(x+n)}f(x+m+n)\\ 
\overline{f(x+p)}f(x+m+p)f(x+n+p)\overline{f(x+m+n+p)} \ .
\end{multline}

These norms were used by Gowers
to control averages along 
arithmetic progressions in a set of positive upper density, 
and adapted by Green and Tao~\cite{GT} for the average along 
arithmetic progressions in the prime numbers.  
Closely related is the ergodic theoretic 
use of seminorms in~\cite{HK} to study 
the ergodic averages, introduced by Furstenberg~\cite{F1}, 
associated with Szemer\'edi's Theorem.  
In~\cite{HK}, we show that the $k$-th seminorm corresponds to a 
$k$-dimensional nilsystem. 
(See~\cite{HK} for the definition and 
precise statement; in the current context, the definition is 
given in  Section~\ref{sec:nilparallel}.)

A natural question is on which sets a parallelepiped 
structure, and thus a Gowers norm, can be defined.  
More interesting is understanding to what extent 
the correspondence with a $k$-step nilpotent 
group can be made in this more general setting.  
We restrict ourselves to the cases $k=2$ and $k=3$ and 
characterize to what extent this correspondence 
can be made precise.  
As the precise definitions of parallelogram and 
parallelepiped structures are postponed until we 
have developed some machinery, we only give a loose 
overview of the results.   
Essentially, the properties included in the definition of a 
parallelepiped structure are exactly those needed in order 
to define a Gowers type norm.

For a two dimensional parallelogram, 
we completely characterize possible parallelogram structures by an 
abelian group (Corollary~\ref{cor:struc-par}).  
This means that a parallelogram structure 
arises from a $1$-step nilpotent group.  For the corresponding three
 dimensional 
case, the situation becomes more complex.  In 
Theorem~\ref{th:structure}, Theorem~\ref{th:split} and 
Corollary~\ref{th:embed} we show under some 
additional hypotheses, a parallelepiped 
structure corresponds to a $2$-step nilpotent group.  However, 
there are examples (Example~\ref{ex:counterex}) for 
which this hypothesis is not satisfied.  In all cases, 
in Theorem~\ref{th:nik-imbedd} we  are able to show that the parallelepiped 
structure can be embedded in a $2$-step nilsystem.  

In Section~\ref{sec:higher}, we outline to what extent these 
results can be carried out in higher dimensions.  

For all sets on which it is possible to define these structures, 
one can naturally define the corresponding Gowers norm $U_k$.   
We expect that these norms should have other applications, outside 
of those already developed by Gowers, Green and Tao, and the authors.

\subsection{Notation}
Parallelogram structures and
parallelepiped structures
are defined as subsets of the Cartesian powers $X^4$ and $X^8$ of 
some sets or groups
and so we introduce some notation. 
 
When $X$ is a set, we write $X\type 2 = X\times X\times X \times X$
 and let 
$X\type 3$ denote the analogous product with $8$ terms. 

A point in  $X\type 2$ 
 is written $\bx=(x_0,x_1,x_2,x_3)$ or 
$\bx = (x_{00}, x_{01}, x_{10}, x_{11})$ and a point in $X\type 3$ is 
written $\ubx=(x_0,\dots,x_7)$ or 
$\ubx=(x_{000},x_{001},\dots,x_{111})$.
More succinctly, we denote $\bx\in X\type 2$ by  
$\bx=(x_i\colon 0\leq i\leq 3)$ or
$\bx=(x_\epsilon\colon \epsilon\in\{0,1\}^2)$,
and use similar notation for points in $X\type 3$.
  
It is convenient to identify  $\{0, 1\}^2$ 
 with the set of vertices of the Euclidean unit 
square. Then 
the second type of notation  allows 
us to view each coordinate of a point $\bx$ of 
$X\type 2$ as lying at the corresponding vertex.

Each Euclidean isometry of the square  permutes the vertices and thus 
the coordinates of $\bx$. The permutations of $X\type 2$ defined in 
this way are called the \emph{Euclidean permutations} of $X\type 2$.
For example, the  maps
$$
\bx\mapsto (x_{10},x_{11}, x_{00}, x_{01}) \text{ and }
\bx\mapsto (x_{10},x_{00}, x_{11}, x_{01})
$$
are Euclidean permutations.
We use the same vocabulary for $X\type 3$, with the Euclidean 
$3$ dimensional unit cube replacing the square.

If $\pi\colon X\to Y$ is a map, by $\pi\type 2\colon X\type 2\to Y\type 2$ 
we mean 
$$\pi\type 2(\bx) = 
\pi(x_{00})\pi(x_{01})\pi(x_{10})\pi(x_{11})\ .$$  
Similarly, $\pi\type 3$ is defined as the corresponding 
map $\colon X\type 3\to Y\type 3$.  

\section{Parallelograms}
\subsection{Definition and examples}
We start with a formal definition of a parallelogram 
structure on an arbitrary set $X$: 

\begin{definition}
\label{def:parallelograms}
Let $X$ be a nonempty set. 
A {\em weak parallelogram structure} on $X$ 
is a subset $\CP$ of $X\type 2$ so that:
\begin{enumerate}
\item
\label{it:def1}
(Equivalence relation) The relation $\sim$ on $X^2$ defined by 
$(x_{00},x_{01})\sim(x_{10},x_{11})$ if and only if 
$(x_{00},x_{01},x_{10},x_{11})\in\CP$
is an equivalence relation. 
\item 
\label{it:def2}
(Symmetry)
If $(x_{00},x_{01},x_{10},x_{11})\in\CP$, 
then $(x_{00},x_{10},x_{01},x_{11})\in\CP$.
\item
\label{it:def3}
(Closing parallelogram property) For all $x_{00},x_{01},x_{10}\in X$, there exists 
$x_{11}\in X$ so that $(x_{00},x_{01},x_{10},x_{11})\in
\CP$.
\end{enumerate}
We say that $\CP$ is a
\emph{strong parallelogram structure} if  the element $x_{11}$ 
in~\ref{it:def3} is unique.
\end{definition}

\begin{definition}
We call the quotient space $X/\sim$ 
the \emph{base} of the structure $\CP$ 
 and denote it by $B$.  
We write the equivalence class 
of an element $(x,y)$ of $X^ 2$ as $\wec{x,y}$.
\end{definition}

The only part of the definition of a parallelogram 
structure that does not appear to be completely natural is the 
transitivity in the equivalence relation $\sim$.  We 
shall justify this assumption later 
(Proposition~\ref{prop:postrans}).  

\begin{lemma}\label{lem:para1}
  Let $X$ be a nonempty set.  
\begin{enumerate}
\item 
\label{rem:reflex1}
For all $x_0,x_1\in X$, $(x_0,x_0,x_1,x_1)\in\CP$.  

\item  
\label{rem:inv-euc} 
$\CP$ is invariant under all Euclidean permutations of 
$X\type 2$.

\item
\label{rem:vectors}
The relation $\sim$ can be rewritten as: 
for $x_{00},x_{01},x_{10}, x_{11}\in X$,
$(x_{00},x_{01},x_{10}, x_{11})\in\CP$ if and only if 
$\wec{x_{00},x_{01}} = \wec{x_{10},x_{11}}$ if and only if
$\wec{x_{00},x_{10}}=\wec{x_{01},x_{11}}$.  

\item 
\label{rem:sameequiv}
All pairs $(x,x)$ with $x\in X$ belong to 
the same $\sim$-equivalence class. 
\end{enumerate}
\end{lemma}

\begin{proof}
\ref{rem:reflex1} Reflexivity of $\sim$ implies that for all $x_0,x_1\in 
X$ we have $(x_0,x_1,x_0,x_1)\in\CP$.  
By property~\ref{it:def2} of the definition,
$(x_0,x_0,x_1,x_1)\in\CP$.
Part~\ref{rem:inv-euc} follows from the symmetry of $\sim$  
and property~\ref{it:def2} of the definition. 
Part~\ref{rem:vectors} follows from the definition of $\sim$ 
and~\ref{rem:inv-euc}, and 
part~\ref{rem:sameequiv} follows immediately from~\ref{rem:reflex1}.  
\end{proof}

\begin{notation}
We denote the equivalence class of all pairs 
$(x,x), x\in B$ by $1_B$.  
\end{notation}

\subsection{Seminorm for parallelograms}
Given a weak parallelogram structure on a set $X$ we can define a seminorm 
similar to the Gowers norm $\norm\cdot_{U_2}$. 

\begin{notation}
If $X$ is a set $\CF(X)$ denotes the space of complex valued functions 
on $X$ with finite support.
\end{notation}
 \begin{proposition}
\label{prop:norme2}
Let $X$ be a nonempty set and let 
$\CP$ be a weak parallelogram structure on $X$. For  any 
function $f\in\CF(X)$ on $X$, we have
\begin{equation}
\label{eq:GN1}
 \sum_{\bx\in\CP}
f(x_{00})\overline{f(x_{01})}\,\overline{f(x_{10})}f(x_{11})\geq 0 \ .
\end{equation}
Letting 
\begin{equation}
\label{eq:GN2}
\norm f_\CP := \Bigl(
\sum_{\bx\in\CP}
f(x_{00})\overline{f(x_{01})}\,\overline{f(x_{10})}f(x_{11})\Bigr)^{1/4}
\ ,
\end{equation}
the map $f\mapsto\norm f_\CP$ is a seminorm on 
$\CF(X)$ and it is a 
norm if and only if the structure $\CP$ is strong.
\end{proposition}

\begin{proof}
We first note that if $F, G$ are functions on $X^2$ with finite 
support, then 
$$
 \sum_{\bx\in\CP}F(x_{00}, x_{01})G(x_{10}, x_{11})=
\sum_{z\in B}\Bigl(\sum_{\substack{(x, y)\in X^2\\ 
\ov{x, y} = z}} F(x, y)\Bigr)\,\Bigl(
\sum_{\substack{(x, y)\in X^2\\ 
\ov{x, y} = z}} G(x, y)\Bigr)\ .
$$
In particular, taking $G=\overline F$,
\begin{equation}
\label{eq:positive}
\sum_{\bx\in\CP}F(x_{00}, x_{01})\overline{F(x_{10}, x_{11})}\geq 0  \ .
\end{equation}
We deduce also that for $F,G\in\CF(X^2)$,
\begin{multline}
\Bigl| \sum_{\bx\in\CP}F(x_{00}, x_{01})G(x_{10}, x_{11})\Bigr|
\\
\leq
\Bigl(\sum_{\bx\in\CP}F(x_{00}, x_{01})\overline{F(x_{10}, x_{11})}
\Bigr)^{1/2}
\Bigl(\sum_{\bx\in\CP}G(x_{00}, x_{01})\overline{G(x_{10}, x_{11})}
\Bigr)^{1/2}
\label{eq:positive2}
\end{multline}
Taking $F(x_{00}, x_{01}) = f(x_{00})\overline{f(x_{01})}$ 
in~\eqref{eq:positive} 
we obtain~\ref{eq:GN2}.  

We now use this to show the {\em Cauchy-Schwarz-Gowers Inequality}: 
for four functions $f_{00}, f_{01}, f_{10}, f_{11}\in
\CF(X)$ 
\begin{equation}
\label{eq:CSG2}
\Bigl\vert\sum_{\bx\in\CP}
f_{00}(x_{00})\overline{f_{01}(x_{01})}\,
\overline{f_{10}(x_{10})}f_{11}(x_{11}) \Bigr\vert\leq
\Vert f_{00}\Vert_{\CP}\,\Vert f_{01}\Vert_{\CP}\,\Vert 
f_{10}\Vert_{\CP}\,
\Vert f_{11}\Vert_{\CP}\ .
\end{equation}
Taking  $F(x_{00}, x_{01}) = f_{00}(x_{00})\overline{f_{01}(x_{01})}$ 
and  $G(x_{10}, x_{11}) = f_{10}(x_{10})\overline{f_{11}(x_{11})}$ 
in~\eqref{eq:positive2} we get that the square of the left hand 
side of~\eqref{eq:CSG2} is bounded by
$$
\sum_{\bx\in\CP}
f_{00}(x_{00})\overline{f_{01}(x_{01})}\,
\overline{f_{00}(x_{01})}f_{01}(x_{11})  \cdot
\sum_{\bx\in\CP}
f_{10}(x_{00})\overline{f_{11}(x_{01})}\,
\overline{f_{10}(x_{10})}f_{11}(x_{11})\ .
$$
By symmetry, the first sum can be rewritten as 
$$
\sum_{\bx\in\CP}
f_{00}(x_{00})
\overline{f_{00}(x_{01})}\,\overline{f_{01}(x_{01})}f_{01}(x_{11})\ .
$$
After a second use of~\eqref{eq:positive2}, we obtain that this sum is bounded by 
$ \norm{f_{00}}_\CP^2\cdot\norm{f_{01}}_\CP^2$. Using the same method 
for the second term we obtain the Cauchy-Schwarz-Gowers Inequality.

Subadditivity of $\Vert\cdot\Vert_{\CP}$ follows easily and 
$\norm\cdot_\CP$ is a seminorm.

Assuming now that $\CP$ is a strong parallelogram structure we show 
that $\norm\cdot_\CP$ is actually a norm. Let $f\in\CF(X)$ be a 
function such that $\norm f_\CP=0$. Let $a$ be an arbitrary point of 
$X$ and $g=\one_{\{a\}}$. By the Cauchy-Schwarz-Gowers Inequality,
$$
0=
 \sum_{\bx\in\CP}g(x_{00})g(x_{01})g(x_{10})f(x_{11})=
\sum_{x_{11}\in X,\;(a,a,a,x_{11})\in\CP}f(x_{11})=f(a)
$$
and so $f$ is identically zero.

Conversely, if the structure is not strong, we claim that 
$\norm\cdot_\CP$ is not a norm. This can be shown directly, but 
in the interest of brevity we postpone the proof until 
Section~\ref{sec:proof-of-norm}, 
after we have developed certain properties of parallelogram 
structures.  
\end{proof}

We justify the assumption of transitivity:
\begin{proposition}
\label{prop:postrans}
Let $X$ be a finite set and let $\CP\subset X^4$ satisfy 
all assumptions of the definition of a strong 
parallelogram structure other than transitivity of $\sim$.  
Assume that the positivity relation~\eqref{eq:positive} is 
satisfied.  Then $\CP$ is a strong parallelogram structure.
\end{proposition}

\begin{proof}  
Assume that $X$ has $n$ elements.  For 
$(x_{0}, x_1)$ and $(x_2, x_3)$ in $X^2$, we define 
$$
M_{(x_0, x_1), (x_2, x_3)} = \begin{cases}
1 & \text{if }(x_0, x_1, x_2, x_3)\in\CP \\
0 & \text{otherwise}  \ .
\end{cases}
$$
This defines a $n^2\times n^2$ matrix $M$.  
Let $\lambda_1, \ldots, \lambda_{n^2}$ be its eigenvalues.  
This matrix is 
symmetric with $1$'s on the diagonal, and 
the sum of the entries for each row is $n$.  
All eigenvalues of $M$ are 
real and by assumption $0\leq \lambda_i \leq n$ 
for $i=1, \ldots, n^2$.  
Then
$$
\sum_{i=1}^{n^2}\lambda_i = \trace(M) = n^2 \text{ and } 
\sum_{i=1}^{n^2}\lambda_i^2 = \trace(M^2) = n^3\ .
$$
Therefore all $\lambda_i$ are either $0$ or $n$, and 
$M^2= n\cdot M$.  Transitivity follows.  
\end{proof}

\subsection{Examples}
We give two examples that illustrate, in a 
sense to be explained, the general behavior of weak and strong 
parallelogram structures:
\begin{example}
\label{ex:one}    
If $G$ is an abelian group (written with multiplicative notation),
 then 
\begin{align*}
 \CP_G &:=\bigl \{\bg=
(g_{00}, g_{01}, g_{10}, g_{11})\in G\type 2\colon 
g_{00}g_{01}^{-1}g_{10}^{-1}g_{11}=1\bigr\}\\
&=\bigl\{(g,gs,gt,gst)\colon g,s,t\in G\bigr\}
\end{align*}
is a strong 
parallelogram structure on $G$.
\end{example}

\begin{example}
\label{ex:two}  Let $X$ be a set, $G$ an abelian group, $\CP_G$ the 
strong parallelogram structure on $G$ defined in Example~\ref{ex:one} and
$\pi\colon X\to G$ a surjection. 
Let $\CP$ be the inverse image of $\CP_G$ under the map $\pi\type 
2\colon X\type 2\to G\type 2$:
\begin{align*}
\CP
&=\{\bx\in X\type 2\colon\pi\type 2\bx\in\CP_G\}\\
&=\bigl\{ \bx\in X\type 2\colon
\pi(x_{00})\,\pi(x_{01})\inv \pi(x_{10})\inv \pi(x_{11})
=1\bigr\}\ .
\end{align*}
Then $\CP$ is a weak parallelogram structure on $X$; it is not strong 
unless $\pi$ is a bijection.
\end{example}

\subsection{Description of parallelogram structures}
\label{subsec:def-paragrammes}
We give a complete description of parallelogram 
structures. 
\begin{lemma}
\label{lem:mult-B}
The set $B$ can be endowed with a multiplication such  that
\begin{equation}
\label{eq:mult-B}
\wec{a,b}\cdot\wec{b,c}=\wec{a,c}\text{ for every }a,b,c\in X\ .
\end{equation}
With this multiplication, $B$ is an abelian group.
\end{lemma}

\begin{proof}
Let $s,t\in B$. Let $a\in X$. 
By part~\ref{it:def3} of the definition of a parallelogram, 
there exists $b\in X$ with 
$\wec{a,b}=s$ and there exists $c\in X$ with 
$\wec{b,c}=t$. We check that $\wec{a,c}$ does not depend on the 
choices of $a,b,c$ but only on $s$ and $t$. 
Let $,a',b',c'\in X$ satisfy $\wec{a',b'}=s$ and
$\wec{c',d'}=t$. By part~\ref{rem:vectors} of Lemma~\ref{lem:para1}, 
$\wec{a,a'}=\wec{b,b'}=\wec{c,c'}$ and thus $\wec{a,c}=\wec{a',c'}$.
It follows immediately that the multiplication in 
$B$ satisfying~\eqref{eq:mult-B} is unique.

By construction, this multiplication is associative and admits the 
 class $1_B$ as unit element. The inverse of the class 
$\wec{a,b}$ is the class $\wec{b,a}$.

We are left with showing that the operation is commutative. 
Given $s,t,a,b,c$  as above, we can choose $d\in X$ such that 
$\wec{a,d}=t$. Then $(a,b,d,c)\in\CP$ and thus 
$\wec{d,c}=\wec{a,b}=s$ and so $t\cdot s=\wec{a,d}\cdot 
\wec{d,c}=\wec{a,c}=s\cdot t$.  
\end{proof}

\begin{notation}
Henceforth we implicitly chose a point $i\in X$ and define a map 
$\pi\colon X\to B$ by
$$
\pi(x)=\wec{i,x}\text{ for every }x\in X\ .
$$
\end{notation}
\begin{corollary}
\label{cor:struc-par}
Any strong parallelogram structure is isomorphic (in the obvious 
sense) to a strong parallelogram structure of the type described in 
Example~\ref{ex:one}. 

Any weak parallelogram structure is isomorphic (in the obvious 
sense) to a weak parallelogram structure of the type described in 
Example~\ref{ex:two}.
\end{corollary}

\begin{proof}
Assume that we have a nonempty set $X$ with parallelogram 
structure $\CP$, quotient space $B$, a point $i\in X$, 
and a map $\pi\colon X\to B$ defined as above.  
For every $\bx\in X\type 2$ we have
\begin{multline*}
\pi(x_{00})\cdot\pi(x_{01})\inv\cdot\pi(x_{10})\inv\cdot\pi(x_{11})=
\wec{i,x_{00}}\cdot\wec{i,x_{01}}\inv\cdot\wec{i,x_{10}}\inv
\cdot\wec{i,x_{11}}\\
=\wec{x_{00},x_{01}}\cdot\wec{x_{10},x_{11}}\inv
\end{multline*}
and this is equal to $1_B$ if and only if $\bx\in\CP$. Therefore the 
structure $\CP$ is equal to the inverse image of $\CP_B$ under 
$\pi\type 2$.

If $\CP$ is strong then $\pi$ is a bijection. Identifying a point of 
$X$ with its image under $\pi$ we get that $\CP$ is defined as in 
Example~\ref{ex:one}.
\end{proof}

The following proposition follows easily from the preceding discussion:
\begin{lemma}\label{lem:para2}
Let $\CP$ be a parallelogram structure on a nonempty set
$X$. 
For $x,y\in X$, the following are equivalent :
\begin{enumerate}
\item $(x,x,x,y)\in\CP$.
\item $\pi(x)=\pi(y)$.
\item $\wec{x,y}=1_B$.
\item For all $a,b,c\in X$ such that $(a,b,c,x)\in\CP$,  we have 
$(a,b,c,y)\in\CP$.
\item There exist $a,b,c\in X$ with  $(a,b,c,x)\in\CP$ and 
$(a,b,c,y)\in\CP$.
\end{enumerate}
\end{lemma}

\begin{notation} We denote the equivalence relation on $X$ defined 
by these conditions 
by $\equiv$ or $\equivp$.
\end{notation}
This relation is  equality if and only if the structure $\CP$ is strong.
Therefore the map $\pi\colon X\to B$ induces a bijection from
 the quotient space of $X/\equivp$ onto $B$. We identify these two 
sets.

\subsection{End of the proof of Proposition~\ref{prop:norme2}}
\label{sec:proof-of-norm}
Assume that the structure $\CP$ is not strong.  
We show that $\norm\cdot_{\CP}$ is not a norm.  
By Lemma~\ref{lem:para2}, there exist 
distinct points $a, b\in X$ with the same image under $\pi$
and furthermore taking any quadruple all of whose 
entries are either $a$ or $b$ is a parallelogram in $\CP$.  
Setting $f=\one_{\{a\}}-\one_{\{b\}}$, the sum in the 
definition of the seminorm has $8$ terms that are $1$ 
and $8$ are $-1$, and so we have a nonzero 
function with $\norm f_\CP = 0$.

\subsection{More examples}
\label{sec:moreexamples}
We continue this section with another example that plays 
a significant role in the sequel. 
To do so, we 
introduce some notation that may seem a bit strange at the moment,
but lends itself easily to the sort of generalization needed later.  

\begin{definition} 
\label{def:G21}
Let $G$ be a group. We write $G\type{2,2}$ for the \emph{diagonal 
subgroup} of $G\type 2$;
$$
G\type{2,2}:=\bigl\{(g,g,g,g)\colon g\in G\bigr\}\ .
$$
 We write $G\type{2,1}$ for the 
subgroup of $G\type 2$ spanned by the elements
$$
(g,g,1,1)\ ;\ (g,1,g,1)\ ; (g,g,g,g)\text{ for }g\in G\ .
$$
$G\type{2,1}$ is called the \emph{two dimensional edge group} of $G$.
\end{definition}

While the set of generators given for the edge group 
is a minimal one, it is not the most natural for understanding 
the name we give the group.  Using the analogy with the Euclidean 
square $\{0,1\}^2$, the set of generators for the edge group consists 
 of all elements of $G\type 2$ where we place $g$'s in entries 
corresponding to an edge of $\{0,1\}^2$ and $1$'s elsewhere.  
This point of view becomes more natural and useful in the generalization 
to three dimensions.

\begin{notation} Let $G$ be a group.
We write $G_2$ for its commutator subgroup. Recall that $G_2$ is the 
subgroup of $G$ spanned by the elements $[g,h]$, $g,h\in G$, where
$[g,h]=ghg\inv h\inv$.
$G_3$ denotes the  second commutator subgroup of $G$, that is, the 
subgroup of $G$ spanned by the elements $[g,u]$ for $g\in G$ and 
$u\in G_2$. 
\end{notation}
By a short computation, we have:
\begin{lemma} \label{lem:G21}
Let $G$ be a group. Then 
\begin{multline}
\label{eq:G21}
G\type{2,1}= \bigl\{\bg\in G\type 2: 
g_{00}g_{01}\inv g_{10}\inv g_{11}\in G_2\bigr\}\\
=\bigl\{ (g,gh,gk,ghku)\colon g,h,k\in G,\ u\in G_2\bigr\}\ .
\end{multline}
In particular, if $G$ is abelian then $G\type{2,1}$ is equal to set 
$\CP_G$ of Example~\ref{ex:one}:
\begin{equation}
\label{eq:G21abel}
G\type{2,1}:=\bigl\{\bg\in G\type 2: 
g_{00}g_{01}\inv g_{10}\inv g_{11}=1\bigr\}
=\{(g,gs,gt,gst)\colon g,s,t\in G\}\ .
\end{equation} 
\end{lemma}

\begin{example} 
\label{ex:three}
Let $G$ be a group  and 
$F$ a subgroup of $G$ containing $G_2$.  
Thus $F$ is normal in 
$G$ and $B=G/F$ is abelian. Let $\pi\colon G\to B$ be the 
natural homomorphism and let $\CP$ be the weak parallelogram structure on $G$ 
defined in Example~\ref{ex:two}:
$$
\CP 
= \bigl\{\bg\in G\type 2\colon \pi\type 2(\bg)\in (G/F)\type{2,1}
= \{\bg\in G\type 2\colon 
g_{00}g_{01}^{-1}g_{10}^{-1}g_{11}\in F\}\ .
$$
It is easy to check that $\CP$ is a subgroup of $G\type 2$ and that
$ \CP= G\type{2, 1}F\type{2}$.  
\end{example}

\section{Parallelepipeds}
\label{sec:parapedes}

\subsection{Notation}
Parallelepipeds are the three dimensional generalization of 
parallelograms, and so naturally arise as subsets of $X\type 3$.  
Recall that we identify  $\{0,1\}^3$ with the set of vertices 
of the unit Euclidean cube.  
Under this  identification, we can naturally associate appropriate 
subsets of  $\{0,1\}^3$ with vertices, edges, or faces of 
the unit cube.  

Thus if $\ubx\in X\type 3$ and $\eta$ is an edge of the unit cube, 
the element
$\bx_\eta :=\{x_\epsilon\colon \epsilon\in\eta\}$ of $X\times X$ is 
called an \emph{edge} of $\ubx$.
Similarly, if $\sigma$ is a  face  of the unit cube, the  element 
$\ubx_\sigma:=\{x_\epsilon\colon \epsilon\in \sigma\}$ of $X\type 3$ 
is called a face of $\ubx$. 
By mapping each vertex of $\sigma$ to a 
vertex of $\{0,1\}^2$ in increasing lexicographic order we can consider 
$\ubx_\sigma$ as an element of $X\type 2$.

In particular, $\bx'= (x_{000}, x_{001}, x_{010}, 
x_{011})$ and $\bx''=(x_{100}, x_{101}, x_{110}, x_{111})$ are opposite 
faces of $\ubx$ and we often write $\bx=(\bx',\bx'')$, naturally identifying 
$X\type 3$ with $X\type 2\times X\type 2$.

\subsection{Definition of a parallelepiped}
\label{subsec:def-parapedes}
\begin{definition}  
\label{def:parallelepiped} 
Assume that $X$ is a nonempty set with a 
weak parallelogram structure $\CP$.  
A {\em weak parallelepiped  structure} $\CQ$ is a 
subset of $X\type 3$ satisfying: 
\begin{enumerate}
\item 
\label{it:faces}
(Parallelograms) 
For every $\bx\in\CQ$ and every face $\sigma$ of $\{0,1\}^3$, $\bx_\sigma \in\CP$.  

\item \label{it:sym} (Symmetries)
$\CQ$ is invariant under all Euclidean permutations of $\{0,1\}^3$.

\item\label{it:equiv}
(Equivalence relation)
The relation $\approx$ (written $\approxq$ in case of ambiguity)
on $\CP$ defined by
$ \bx'\approx\bx''$ if and only if $(\bx',\bx'')\in\CQ$
is an equivalence relation.

\item 
\label{it:closing}(Closing parallelepiped property)
If $x_{000}, x_{001}, x_{010}, x_{011}, x_{100}, x_{101}, x_{110}$ are
seven points of $X$ satisfying 
$$
(x_{000}, x_{001}, x_{010}, x_{011}), 
(x_{000},  x_{010},  x_{100}, x_{110}), \text{ and  }
(x_{000}, x_{001}, x_{100}, x_{101})\in\CP \ ,
$$
then there exists $x_{111}\in X$ such that 
$$
(x_{000}, x_{001}, x_{010}, x_{011}, x_{100}, x_{101}, 
x_{110},x_{111})\in\CQ\ .
$$
\end{enumerate}

We say that $\CQ$ is a {\em strong parallelepiped structure} if the element 
$x_{111}$  in~\ref{it:closing} is unique.  

A parallelepiped structure on $X$ with parallelograms $\CP$ and 
parallelepipeds $\CQ$ is denoted by $(\CP, \CQ)$.  
\end{definition}

\begin{notation}
We denote the quotient space $\CP/\approx$ by $P$ and denote 
the equivalence class of a parallelogram $\bx\in\CP$ by 
$\br\bx$ or $\br\bx_\CQ$.
\end{notation}

We start with some properties that follow immediately from the definition: 

Let $x_{000},\dots,x_{110}$ be seven points in $X$ 
satisfying the hypotheses of condition~\ref{it:closing}. Then there 
exists $x_{111}\in X$ such that $(x_{010},x_{110}, x_{011}, 
x_{111})\in\CP$ and one can easily check that all faces of 
$\ubx=(x_{000}, x_{001}, x_{010}, x_{011}, x_{100}, x_{101}, 
x_{110},x_{111})$ belong to $\CP$. Therefore 
condition~\ref{it:closing} can be rewritten as :
\begin{enumerate}
\setcounter{enumi}{4}
\item\label{rem:closing2}
Let $\ubx\in X\type 3$ be such that each of its faces belongs to $\CP$. 
Then there exists $x'_{111}\in X$ such that
$(x_{000}, x_{001}, x_{010}, x_{011}, x_{100}, x_{101}, 
x_{110},x'_{111})\in\CQ$.
\end{enumerate}

We also note that 
some of the conditions in the definition of 
a parallelepiped are redundant: if $\CQ$ 
is invariant under the permutations given in condition~\ref{it:sym} and if the 
relation $\approx$ is  transitive,  then $\approx$ 
is an equivalence relation.

\begin{lemma}
\label{prop:trivial}
Let $X$ be a nonempty set with parallelepiped structure $(\CP, \CQ)$.  Then 
\begin{enumerate}
\item 
\label{eq:xyxy}
For $x,y\in X$, $\br{x,y,x,y}$ depends only on 
$\wec{x,y}$.  

\item 
\label{prop:sameequiv}
All parallelograms of the form 
$(a,a,a,a)$ for some  $a\in X$ belong to the same $\approx$ equivalence class.
\end{enumerate}
\end{lemma}

\begin{proof}
If $\ov{x,y}=\ov{x',y'}$, then $\bx=(x,y,x',y')\in \CP$. 
By reflexivity of the relation $\approx$, we have $(\bx,\bx)\in\CQ$.
 By Part~\ref{it:sym} of the definition of a parallelepiped, 
$(x,y,x,y,x',y',x',y')\in\CQ$ and $\br{x,y,x,y}=\br{x',y',x',y'}$.
  Part~\ref{prop:sameequiv} is an immediate corollary.  
\end{proof}

\begin{notation}
The common equivalence class in Part~\ref{prop:sameequiv} 
is called the \emph{trivial class} and is denoted by 
$1_P$.
\end{notation}

We note that we did not make 
use of the closing parallelepiped property, part~\ref{it:closing} 
of the definition of a parallelepiped property in this proposition.

\subsection{Reduction to strong parallelepiped structures}
For the moment we assume that $\CQ$ satisfies the first three 
properties of the definition, but do not assume the closing 
parallelepiped property (Part~\ref{it:closing} of the definition).

\begin{proposition}
\label{prop:trivialequiv}
Let $X$ be a nonempty set with parallelogram structure $\CP$ and a 
subset $Q$ of $X\type 3$
only satisfying 
properties~\ref{it:faces}, \ref{it:sym} and~\ref{it:equiv}
of definition~\ref{def:parallelepiped} of a parallelepiped 
structure.  
For $x,x'\in X$ the following are equivalent:
\begin{enumerate}
\item $(x,x,x,x')\in\CP$ and $\br{x,x,x,x'}=1_P$.
\item For all $a,b,c\in X$ such that $(a,b,c,x)\in\CP$, we have 
$(a,b,c,x')\in\CP$ and
$\br{a,b,c,x}=\br{a,b,c,x'}$.
\item
There exist $a,b,c\in X$ with $(a,b,c,x)\in\CP$, $(a,b,c,x')\in\CP$ 
 and
$\br{a,b,c,x}=\br{a,b,c,x'}$.
\end{enumerate}
\end{proposition}

\begin{proof}Each of the three 
properties implies that $\ov{x,x'}=1_B$ and thus that
 $(x,x,x,x')\in\CP$. 
Assume first that $\br{x,x,x,x'}= 1_P$. 
Let $a,b,c\in X$ be such that $(a,b,c,x)\in\CP$.
Since $\br{c,c,c,c}=1_P=\br{x,x,x,x'}$,  by symmetry  we have that 
$\br{c,x,c,x}=\br{c,x,c,x'}$. 
By Proposition~\ref{prop:trivial}, part~\ref{eq:xyxy} we have
$\br{a,b,a,b}=\br{c,x,c,x}=\br{c,x,c,x'}$ and thus 
$\br{a,b,c,x}=\br{a,b,c,x'}$ and $x,x'$ satisfy the second property.  

The second condition trivially implies the third.  

Assume now that there exist $a,b,c\in X$ with  $(a,b,c,x)\in\CP$ and
$\br{a,b,c,x}=\br{a,b,c,x'}$. The same argument implies that 
$\br{x,x,x,x'}=\br{x,x,x,x}=1_P$, and our claim is proved.
\end{proof}

\begin{notation}
For $x,x'\in X$, we write $x\equivq x'$ if $x, x'$ satisfy 
any of the three equivalent properties in Proposition~\ref{prop:trivialequiv}. 
\end{notation}
 
The third property implies that $\equivq$ is an 
equivalence relation. Moreover, $\CQ$ is saturated for this relation 
meaning that for $\ubx,\ubx'\in X\type 3$ we have
if $\ubx\in\CQ$ and $x_\epsilon\equivq x'_\epsilon$ for every 
$\epsilon\in\{0,1\}^3$, then $\ubx'\in\CQ$.  
In particular, the structure $(\CP,\CQ)$ is strong if and only if 
relation $\equivq$ is the equality.

Let $Y$ be the quotient space $X/\equivq$ and $r\colon X\to Y$ the 
quotient map. Let $\CP_Y$ and $\CQ_Y$ be the images of $\CP$ and 
$\CQ$ under $r\type 2$ and $r\type 3$, respectively. 
Then $\CP$ and $\CQ$ are the inverse images of $\CP_Y$ and $\CQ_Y$ 
under these maps. 

If $\CQ$ satisfies the closing 
parallelepiped property, then $\CQ_Y$ also satisfies this property and 
$(\CP_Y,\CQ_Y)$ is a parallelepiped structure on $Y$.  This 
structure is  strong because clearly the relation $\equivqy$ is the 
identity.

We have thus shown:
\begin{proposition}\label{prop:weakstrong}
Every parallelepiped structure on a nonempty set is the inverse 
image of a strong parallelepiped structure.
\end{proposition}

Thus the study of parallelepiped structures 
reduces to the study of strong ones and so 
in  the sequel, we generally consider strong parallelepiped 
structure.

\subsection{The seminorm associated to a a parallelepiped structure}
Given a parallelepiped structure, one can define  a third 
(Gowers) norm, giving us a three dimensional version 
of Proposition~\ref{prop:norme2}.  We include a proof 
for completeness, but note that we have included in 
the definition of a parallelepiped structure all the properties needed 
to parallel the steps in Gowers's original proof.  
(The notation is given in the introduction.)

\begin{proposition}
\label{prop:norme3}
Let $(\CP,\CQ)$ be a parallelepiped structure on the set $X$.
For every $f\in\CF(X)$,
\begin{equation}
\label{eq:gowers1}
 \sum_{\ubx\in\CQ}\prod_{\epsilon\in\{0,1\}^3}C^{|\epsilon|}f(x_\epsilon)\geq 0
\end{equation}
and thus we can define
\begin{equation}
 \norm f_\CQ :=
\Bigl( \sum_{\ubx\in\CQ}\prod_{\epsilon\in\{0,1\}^3}
C^{|\epsilon|}f(x_\epsilon)\Bigr)^{1/8}\ .
\end{equation}
The map $f\mapsto \norm f_\CQ$ is a seminorm on $\CF(X)$ and it is a 
norm if and only if the structure $(\CP,\CQ)$ is strong.
\end{proposition}

We omit the proof, as it is exactly the same as the proof 
of the two dimensional version, Proposition~\ref{prop:norme2}.  
For the converse implication to show that if 
we have a norm the structure is strong, we use 
Proposition~\ref{prop:trivialequiv}
instead of Lemma~\ref{lem:para2}.

Similarly, the assumption of transitivity of the relation 
$\approx$ in the definition of parallelepipeds is related 
to positivity, as for parallelograms in 
Proposition~\ref{prop:postrans}:
\begin{proposition}
\label{prop:postrans3}
Let $X$ be a finite set, let $\CP$ be a parallelogram 
structure on $X$ and let $\CQ\subset X\type 3$ satisfy 
all assumptions of the definition of a strong 
parallelepiped structure other than transitivity of $\approx$.  
Assume that the following positivity relation holds: for 
every function $F$ on $X^4$ with finite support,  
$$
\sum_{\ubx\in\CQ} f(\bx')\overline{f(\bx'')} \geq 0 \ .
$$
Then $\CQ$ is a strong parallelogram structure.
\end{proposition}
We omit the proof.  

\subsection{First examples: Abelian parallelepiped structures}
\label{subsec:abelian}

We begin with some examples of parallelepiped structure on abelian 
groups. We need some  notation and 
definitions similar to the definitions introduced for
Definition~\ref{def:G21} in the two-dimensional case.

\begin{notation}
Let $G$ be a group.  
For a nonempty subset $\alpha$ of $\{0,1\}^3$ and $g\in G$, 
we write $g\type{3,\alpha}$ 
for the element of $G\type 3$ given by
$$
 \bigl(g\type{3,\alpha}\bigr)_\epsilon=\begin{cases}
g & \text{if }\epsilon\in \alpha \\
1 & \text{otherwise\ .}
\end{cases}
$$

Note that the elements $g\type {3,\eta}$, 
$g\in G$, $\eta\in\{0,1\}^3$, generate $G\type 3$.
If $\alpha=\{0,1\}^3$ then $g\type {3,\alpha}$ is the diagonal 
element $(g,g,\dots,g)$.  
\end{notation}

\begin{definition} Let $G$ be a group.

The \emph{diagonal group}  $G\type{3,3}$ is the subgroup of $G\type 
3$ consisting in elements of the form  $(g,g,\dots,g)$ for $g\in G$.

The \emph{edge group} $G\type{3,1}$  is the subgroup 
of $G\type 3$ spanned by the elements of the form $g\type{3,e}$ where 
$g\in G$ and $e$ is an edge of the cube $\{0,1\}^3$.

The \emph{face group} $G\type{3,2}$ is the subgroup of $G\type 3$ 
spanned by the elements of the form $g\type{3,f}$ where $g\in G$ and 
$f$ is a face of the cube $\{0,1\}^3$.
\end{definition}

The following Proposition follows immediately:.
\begin{proposition}
If $G$ is an abelian group, then 
\begin{equation}
\label{eq:G31abelian}
 G\type{3,1}=
\Bigl\{\ubg\in G\type 3\colon 
\prod_{\epsilon\in\{0,1\}^3}g_\epsilon^{(-1)^{|\epsilon|}}=1
\Bigr\}
\end{equation}
and 
$$
 G\type{3,2}=\bigl\{\ubg\in G\type 3\colon\text{ every face of $\ubg$ belongs 
to }G\type{2,1}\bigr\}\ .
$$
\end{proposition}

\begin{example}
\label{ex:abel1}
Let $G$ be an abelian group, let $\CP=G\type{2,1}$ be the strong 
parallelogram structure 
defined as in Example~\ref{ex:one}  
and let $\CQ=G\type{3,2}$.  
Then $(\CP,\CQ)$ is a (strong) parallelepiped structure on 
$G$.  
\end{example}

\begin{example}\label{ex:abel2}
Let $G$ be an abelian group and $F$ a subgroup of $G$. Define
$$
\CP=G\type{2,1}F\type 2\text{ and }\CQ=G\type{3,2}F\type{3,1}\ .
$$
Then $(\CP,\CQ)$ is a (strong) parallelepiped structure on $X$.
\end{example}
This assertion is a particular case of a more general statement 
(Proposition~\ref{ex:group2}) and so we omit the proof.

\subsection{Some nonabelian examples}
\label{sec:nonabelian}
More interesting are  parallelepiped structures on nonabelian groups.  
It is here that the structures take on nontrivial properties.
 
We begin with an elementary remark:
\begin{remark}\label{rem:commut}
Let $G$ be a group and let $\alpha,\beta$ be two subsets of the cube 
$\{0,1\}^3$. Then, for every $g,h\in G$, the commutator of the two 
elements $g\type{3,\alpha}$ and $h\type{3,\beta}$ of $G\type 3$ is
$$
 \bigl[g\type{3,\alpha},h\type{3,\beta}\bigr]=[g,h]\type{3,\alpha\cap\beta}\ .
$$
\end{remark}
\begin{lemma}
\label{lem:comut}
Let $G$ be a  group. Then:
\begin{enumerate}
\item\label{it:comut1}
 $G\type{3,2}\supset G_2\type{3,1}\supset G_3\type 3$.

\item\label{it:comut3}
Let $g\in G$ and $\eta\in\{0,1\}^3$. Then $g\type{3,\eta}\in G\type 
{3,2}$ if and only if $g\in G_3$.

\item\label{it:comut5} 
$G\type{2,2}G_2\type{2,1}G_3\type 2$ is a normal 
subgroup of $G\type{2,1}$. 

Moreover, under the identification $G\type 3=G\type 2\times G\type 2$,
\begin{equation}
\label{eq:gprime}
 G\type{3,2}=\bigl\{ \ubg=(\bg',\bg'')\in G\type {2,1}\times G\type 
{2,1}\colon g'{g''}\inv\in 
G\type{2,1}G_2\type{2,1}G_3\type 2\ .
\end{equation}
\end{enumerate}
\end{lemma}

In a more general context, the proof is 
contained in  sections $5$ and $11$ of~\cite{HK}.  
The idea is to find the natural setting in which these 
cubic structures form a group, much as Hall~\cite{H} and 
Petresco~\cite{P} (see also~\cite{Laz}), 
and later Leibman~\cite{L}, did for arithmetic progressions.
These groups also arise in~\cite{GT3}, and Green and Tao 
refer to $G\type{3,2}$ as the Hall-Petresco cube group.
For completeness, we summarize the argument.  

\begin{proof}  
\noindent\emph{Part~\ref{it:comut1}}
Each face of $\{0,1\}^3$ is the union of two edges. 
If $\sigma$ and $\tau$ are two  faces of $\{0,1\}^3$, then 
$\sigma\cap\tau$ is a face, an edge, or the empty set.
Conversely, each edge can be written as 
the intersection of two faces. 
By Remark~\ref{rem:commut}, 
the commutator subgroup of $G\type{3,2}$ is
therefore equal to
$G_2\type{3,1}$. By 
a similar argument, the second commutator subgroup of $G\type{3,2}$ is 
$G_3\type 3$, and Part~\ref{it:comut1} follows.

\noindent\emph{Part~\ref{it:comut5}}
We write  $\one=(1,1,1,1)\in G\type 2$ and 
$K=G\type{2,2}G_2\type{2,1}G_3\type 2$.
One can check directly that $K$ is a normal subgroup of 
$G\type{2,1}$. 

It follows immediately from the definition of $G\type{3,2}$ that for 
$\ubg=(\bg',\bg'')\in G\type{3,2}$ we have that $\bg'$ and $\bg''$ 
belong to $G\type{2,1}$ and that for $\bg\in G\type{2,1}$ we have 
$(\bg,\bg)\in G\type{3,2}$. Therefore 
$$
 G\type{3,2}=
\bigl\{(\bg',\bg'')\in G\type{2,1}\times G\type{2,1} 
\colon g'{g''}\inv\in L\bigr\}\ ,
$$
where 
$$
 L=\bigl\{\bg\in G\type{2,1}\colon (\bg,\one)\in G\type{3,2}\}\ .
$$
We are left with checking that $L=K$.  
The inclusion $K\subset L$ follows 
 immediately from the definition of $G\type{3,2}$ and 
Part~\ref{it:comut1}. Moreover, the subset of $G\type 3$ on the right 
hand side of~\eqref{eq:gprime}
is a subgroup of $G\type 3$ containing the generators of 
$G\type{3,2}$ and 
thus containing this group. This implies that $K\supset 
L$.

\noindent\emph{Part~\ref{it:comut3}} If $g\in G_3$ and 
$\eta\in\{0,1\}^3$, then $g\type{3,\eta}\in G\type{3,2}$ by 
Part~\ref{it:comut1}. 
Conversely, let $g\in G$ and assume that 
$g\type{3,\eta}\in G\type{3,2}$ for some vertex $\eta$. We have to 
show that $g\in G_3$. By the symmetries of $G\type{3,2}$, we 
can restrict to the case that $\eta=111$.

By Part~\ref{it:comut5}, $(1,1,1,g)$ belongs to
$G\type{2,2}G_2\type{2,1}G_3\type 2$ and we can write
$$
 (1,1,1,g)=(h,h,h,h).\bu.\bv\text{ with }h\in G,\ 
\bu\in G_2\type {2,1}\text{ and }\bv\in G_3\type 3\ .
$$
Looking at the congruences modulo $G_2$ of the first coordinate, 
we have that $h\in G_2$ and thus $(h,h,h,h)\in G_2\type{2,1}$.
Substituting $(h,h,h,h).\bu$ for $\bu$, we reduce to the case that 
$h=1$. 

Recall that $G_3$ is a normal subgroup of $G_2$ and that $G_2/G_3$ is 
abelian. Let $\bar\bu$ be the element of $(G_2/G_3)\type 2$ obtained by reducing each 
coordinate of $\bu$ modulo $G_3$. Then $\bar\bu$ belongs to 
$(G_2/G_3)\type{2,1}$, its first three coordinates are equal to $1$ 
and by Lemma~\ref{lem:G21} its last coordinate is also equal to $1$.
This means that $u_{11}\in G_3$ and it follows that $g\in G_3$.
\end{proof}

We now turn to several nonabelian generalizations of the previous 
examples. The first one generalizes
Example~\ref{ex:abel1}.

\begin{proposition}
\label{ex:group1}
Let $G$ be a  group, $\CP=G\type{2,1}$, and $\CQ=G\type{3,2}$. Then
$(\CP,\CQ)$ is a  parallelepiped structure on $G$ and this structure 
is strong if and only if $G$ is $2$-step nilpotent.
\end{proposition}

We postpone the proof until after a second example of a group 
parallelepiped structure, which generalizes 
Example~\ref{ex:abel2} to the nonabelian setting.

\begin{proposition}
\label{ex:group2}
Let $G$ be a group, $F$ a subgroup of $G$  with
\begin{equation}
\label{eq:hypgroup2}
 G_2\subset F\subset Z(G)\ .
\end{equation}
We define (as in Proposition~\ref{ex:three})
$\CP=G\type{2,1}F\type 2 $ 
and $\CQ=G\type{3,2}F\type{3,1}$.  
Then $(\CP,\CQ)$ is a strong parallelepiped structure on $G$.
\end{proposition}

Note that condition~\eqref{eq:hypgroup2} implies in particular that $G$ is a 
$2$-step nilpotent group.

\begin{proof}
The symmetries of $\CQ$ are obvious. 
$\CQ$ is clearly a subgroup of $G\type 3$ and it is not 
difficult 
to deduce form part~\ref{it:comut5} of 
Lemma~\ref{lem:comut}  that $G\type{2,2}F\type{2,1}$ is normal in 
$G\type{2,1}F\type 2$ and that
\begin{equation}
\label{eq:CQ2}
 \CQ=\{(\bg',\bg'')\in G\type 2\times G\type 2\colon
\bg'\in G\type{2,1}F\type 2,\quad \bg'{\bg''}\inv\in 
G\type{2,2}F\type{2,1}\}\ .
\end{equation}
It follows that for every element $\ubg = (\bg', \bg'')\in\CQ$, 
the ``first'' face $\bg'$ of 
$\ubg$ belongs to $\CP$.  Thus by symmetry, all 
faces belong to $\CP$.  It also follows 
from~\eqref{eq:CQ2} that $\approxq$ is an equivalence relation on 
$\CP$.

Let $\CQ'$ be the subgroup of $G\type 3$ consisting of elements $\ubg$ 
such that each face of $\ubg$ belongs to $\CP$. Then $\CQ'$ contains 
$F\type 3$ and the quotient group $\CQ'/F\type 3$ is the abelian 
parallelepiped  structure $(G/F)\type{3,2}$ on $G/F$ defined in 
Example~\ref{ex:abel1} and thus $\CQ'=G\type{3,2}F\type 3$.

Let $\ubg\in\CQ'$. We write $\ubg=\ubj\ubh$ with 
$\ubj\in G\type{3,1}$ and $\ubh$ of $F\type 3$. 
Since $F$ is abelian, by~\eqref{eq:G31abelian} 
 there exists an element 
$\ubh'\in F\type{3,1}$ that is equal to $\ubh$ other than 
in the last coordinate.
Thus $\ubg':=\ubj \ubh'\in\CQ$ and coincides with $\ubg$ other than 
in the last coordinate.  
Therefore, condition~\ref{rem:closing2} is satisfied.

Finally we show that the structure is strong.   Assume not.  
Since $\CQ$ is a group, there 
exists $\ubg\in\CQ$ of the form $(1,1,\dots, 1,g)$ for some $g\in G$ 
not equal to $1$. 
By~\eqref{eq:CQ2}, $(1,1,1,g)\in G\type{2,2}F\type{2,1}$. 
Applying formula~\eqref{eq:G21abel} for $F\type{2,1}$,  we have that 
$g=1$, a contradiction.
 \end{proof}

We now return to 
Proposition~\ref{ex:group1}, and show that $(\CP,\CQ)$ as defined in this 
example is a parallelepiped structure on $G$.

\begin{proof}[Proof of Proposition~\ref{ex:group1}]  
First consider the case that $G$ is  $2$-step nilpotent.
Then the hypotheses of Proposition~\ref{ex:group2} are satisfied with 
$F=G_2$. In order to show that $(\CP,\CQ)$ is a strong parallelepiped 
structure we  check that $\CP$ and $\CQ$ are equal to the sets defined 
in this Proposition.
 By 
Lemma~\ref{lem:G21},  we have $ G_2\type 2\subset G\type{2,1}$ and thus 
$G\type{2,1}G_2\type 2=G\type{2,1}=\CP$; by part~\ref{it:comut1} of 
Lemma~\ref{lem:comut} we have that 
$G\type{3,2}G_2\type{3,1}=G\type{3,2}=\CQ$. 

We now turn to the general case. Recall that $G/G_3$ is a $2$-step 
nilpotent group.  Since the group $\CP$ contains 
$G_3\type 3$ and $\CQ$ contains $G_3\type 3$, 
$\CP$ is the inverse image of $(G/G_3)\type{2,1}$  in 
$G\type 2$ and $\CQ$ is the inverse image of 
$(G/G_3)\type{3,2}$ in $G\type 3$.  The assertion follows 
immediately.
\end{proof}

\subsection{Nilparallelepiped structures}
\label{sec:nilparallel}

Let $G$ be a group and $F$ a subgroup of $G$ satisfying 
$G_2\subset F\subset Z(G)$ (hypotheses~\eqref{eq:hypgroup2}
of Proposition~\ref{ex:group2}).  
Let  $\Gamma$ be a subgroup, not necessarily normal,  of $G$.  
We build a structure on the coset space $X= G/\Gamma$.
By substituting $G/(\Gamma\cap F)$ for $G$, $F/(\Gamma\cap F)$ for 
$F$, and $\Gamma/(\Gamma\cap F)$ for $\Gamma$ we reduce to the case 
that $\Gamma\cap F=\{1\}$.

\begin{proposition}
\label{prop:nil}
Let $G$ be a group and $F$ a subgroup of $G$ satisfying 
hypothesis~\eqref{eq:hypgroup2}
of Proposition~\ref{ex:group2}:
\begin{equation*}
\tag{\ref{eq:hypgroup2}}
G_2\subset F\subset Z(G)
\end{equation*}
and let $\Gamma$ be a subgroup of $G$ 
with
\begin{equation}
\label{eq:hypnil}
\Gamma\cap F=\{1\}\ .
\end{equation}
Let 
$\CP=G\type{2,1}F\type 2$, $\CQ=G\type{3,2}F\type{3,1}$, 
$\pi\colon G\to G/\Gamma$ be projection, and let 
$$
 \CP_X=\pi\type 2(\CP)\ \ \text{ and } \  \CQ_X=\pi\type 3(\CQ)\ .
$$
Then $(\CP_X,\CQ_X)$ is a strong parallelepiped structure on $X$.
\end{proposition}

A parallelepiped structure defined by $G, F, \Gamma, X, \CP_X, \CQ_X$ 
as in Proposition~\ref{prop:nil} is called a 
\emph{nilmanifold parallelepiped structure}, and more succinctly we 
refer to it as a nilparallelepiped structure.

Before the proof of Proposition~\ref{prop:nil}, we have a lemma.
\begin{lemma}  
\label{cl:closing3} 
Maintaining notation as in Proposition~\ref{prop:nil}, if 
$\bg=(g_{00},g_{01},g_{10},g_{11})\in G\type 2$ is
such that $\pi\type 2(\bg)\in\CP_X$, then there exists $g'_{11}\in G$ 
such that $\pi(g'_{11})=\pi(g_{11})$ and 
$(g_{00},g_{01},g_{10},g'_{11})\in\CP$.
\end{lemma}

\begin{proof}
By definition there exists $\bh\in \CP$ with $\pi\type 
2(\bh)=\pi\type 2(\bg)$ and thus
there exists $\bgamma\in \Gamma\type 2$ with $\bg=\bh\bgamma$.
As $\Gamma$ is abelian  it follows from 
Lemma~\ref{lem:G21} that there exists $\theta\in\Gamma$ with 
$(\gamma_{00},\gamma_{01},\gamma_{10},\gamma_{11}\theta)\in 
\Gamma\type{2,1}$ and the point $g'_{11}=g_{11}\theta$ satisfies the 
announced properties.  
\end{proof}

\begin{proof}[Proof of Proposition~\ref{prop:nil}]
$\CP_X$ is clearly the weak parallelogram structure defined as in 
Example~\ref{ex:two} by the natural projection of $X$ on the base 
group $B=G/F\Gamma$. 
If $\ubx\in\CQ_X$ then every face of $\ubx$ belongs to $\CP_X$ by 
definition. The symmetries of $\CQ$ are obvious.

Let $\bx,\by,\bz$ be three parallelograms in $\CP_X$, such that 
$(\bx,\by)$ and $(\by,\bz)$ belong to $\CQ_X$.  This means that there 
exist two parallelepipeds $\ubg=(\bg',\bg'')$ and $\ubh=(\bh',\bh'')$ in 
$\CQ$ such that 
$$
 \pi\type 2(\bg')=\bx\ ;\ \pi\type 2(\bg'')=\pi\type 2(\bh')=\by\text{ and }
\pi\type 2(\bh'')=\bz\ .
$$
 Let $\bgamma={\bg''}\inv{\bh'}$. By~\eqref{eq:CQ2},
we have $\bgamma\in\Gamma\type 2\cap G\type{2,1}F\type 2$ and thus
$\gamma_{00}\gamma_{01}\inv \gamma_{10}\inv \gamma_{11}\in \Gamma\cap F=\{1\}$ 
and so $\bgamma\in \Gamma\type{2,1}$. By part~\ref{it:comut5} of 
Lemma~\ref{lem:comut} again, $(\bg'\bgamma,\bh')\in\CQ$ and thus 
$(\bg'\bgamma,\bh'')\in\CQ$ by transitivity. The projection of this 
parallelepiped on $X$ 
is $(\bx,\bz)$ and thus $(\bx,\bz)\in\CQ$. This shows that the 
relation $\approx_{\CQ_X}$ on $\CP_X$ is transitive and we deduce 
that it is an equivalence relation.

Consider now a point $\ubx\in X\type 3$ such that each of its faces 
belongs to $\CP_X$. By Lemma~\ref{cl:closing3},  there exists 
$\ubg\in G\type 3$ with $\pi\type 3(\ubg)=\ubx$ and such that
 $\ubg_f\in\CP$ when $f$ is any of the three
faces of $\{0,1\}^3$ containing $000$ and also for the face 
$\{100,101,110,111\}$. It follows that the two remaining faces of 
$\ubg$ also are parallelograms
of $G$.  Since $\CQ$ is a parallelepiped structure on $G$, we can modify 
(again using Lemma~\ref{cl:closing3}) 
$g_{111}$ in order to get a parallelepiped in $G$.  
Projecting this parallelepiped 
on $X$, we get a parallelepiped in $\CQ_X$ coinciding with $\ubx$ other than 
in the last coordinate.

Thus we have that $\CQ_X$ is a parallelepiped structure on $X$ 
and we are left with showing it is strong. If not, 
there exist two parallelepipeds $\ubx,\uby\in\CQ_X$ with 
$x_\epsilon=y_\epsilon$ for every $\epsilon\neq 111$ and $x_{111}\neq 
y_{111}$. Let $\ubg,\ubh$ be two parallelepipeds in $G$ with 
$\pi\type 3(\ubg)=\ubx$ and $\pi\type 3(\ubh)=\uby$. Writing 
$\ubu=\ubg\inv\ubh$ we have $ \ubu\in\CQ$ with $u_\epsilon\in\Gamma$ 
for every $\epsilon\neq 111$ and $u_{111}\notin\Gamma$.
By part~\ref{it:comut5} of Lemma~\ref{lem:comut}, there exists $\bv\in 
G\type{2,2}G_2\type{2,1}$ with $v_\epsilon\in\Gamma$ for $\epsilon\neq 
11$ and $v_{11}\notin\Gamma$. By Lemma~\ref{lem:G21}, $\bv$ can be 
written as $\bv=(v,vs,vt,vst)$ with $v\in G$ and $s,t\in G_2$ and we get 
that $v,vs,vt$ and thus $vst$ belong to $\Gamma$, a contradiction.
\end{proof}

\subsection{Nilpotent groups appear}

All the examples of strong parallelepiped structures
considered thus far have a striking feature in common. 
There exist a $2$ step nilpotent group $G$ acting transitively on $X$.
We explore this further.  

The group $G\type 2$ acts on $X\type 2$ in a 
natural way and we write $(\bg,\bx)\mapsto\bg\cdot\bx$ for this action.
In all examples thus given 
it can be checked that for every parallelogram 
$\bx\in\CP$ and every $g\in G$ we have that $g\type 2\cdot\bx$ is a 
parallelogram, equivalent to $\CP$ under the relation $\approx$.

This is not merely coincidence; it reflects an underlying 
structure, so long as certain algebraic considerations are avoided.  
We denote the action of $G$ on $X$ by
$(g,x)\mapsto g\cdot x$. 
 This 
motivates the following definition.  

We need some notation. Let $g\colon x\mapsto g\cdot x$ be a 
transformation of $X$. We recall that $g\type 2$ and $g\type 3$ are 
the diagonal transformations of $X\type 2$ and $X\type 3$, 
respectively. More generally, if $\alpha$ is a subset of the cube 
$\{0,1\}^3$, $g\type {3,\alpha}$ is the transformation of $X\type 3$ 
given by
$$
\Bigl( g\type{3,\alpha}\cdot\ubx\Bigr)_\epsilon=
\begin{cases} g\cdot x_\epsilon & \text{if }\epsilon\in\alpha\\
x_\epsilon &\text{otherwise.}
\end{cases}
$$
This notation is coherent with the notation introduced in 
Section~\ref{subsec:abelian}.

\begin{definition}
\label{def:struc}
Let $(\CP,\CQ)$ be a strong parallelepiped structure on a nonempty set $X$. 
The \emph{structure group} of $(\CP,\CQ)$, written $\CG$ or $\CG_Q$, 
is the group of bijections $x\mapsto g\cdot x$ of $X$ such that
for every parallelogram $\bx\in\CP$, $g\type 2\cdot\bx$
is a parallelogram  and $g\type 2\cdot\bx\approx\bx$.  
\end{definition}

We can rephrase this condition:

\begin{proposition}  If $(\CP, \CQ)$ is a strong parallelepiped 
structure on a nonempty set $X$, then the 
structure group $\CG$ is the group of 
bijections $x\mapsto g\cdot x$ of $X$ such that 
for every parallelepiped $\ubx\in\CQ$ and every face $f$ of 
$\{0,1\}^3$, we have $g\type {3,f}\cdot \ubx\in\CQ$.  
\end{proposition}

\begin{proof} Let $\phi$ be the face 
$\{\epsilon\in\{0,1\}^3\colon\epsilon_3=1\}$.

Assume first that $g\in\CG$ and let $\ubx\in\CQ$.  As usual, we write 
$\ubx=(\bx',\bx'')$ with $\bx',\bx''\in\CP$ and
we have $g\type{3,\phi}\cdot \ubx=(\bx',g\type 2\cdot\bx'')$.
By hypothesis, 
$g\type 2\cdot \bx''\in\CP$ and  $g\type 2\cdot \bx''\approx \bx''$, 
meaning that 
$(\bx'',g\type 2\cdot \bx'')\in\CQ$. By transitivity,
$g\type{3,\phi}\cdot\ubx\in\CQ$.
By symmetry, the same result also holds for the other faces of the 
cube.

Assume now that $x\mapsto g\cdot x$ is a bijection of $X$ such that 
for every $\ubx\in\CQ$ and every face $f$ we have 
$g\type{3,f}\cdot\ubx\in\CQ$. 
Let $\bx\in\CP$. Then $\ubx:=(\bx,\bx)\in\CQ$ and so $(\bx,g\type 
2\cdot\bx)=g\type{3,\phi}\cdot\ubx\in\CQ$.  Thus $g\type 
2\cdot\bx\in\CP$ and $g\type 2\cdot\bx\approx\bx$.
\end{proof}

\begin{proposition}
\label{prop:nil2}
Let  $(\CP,\CQ)$ be a strong parallelepiped structure on a nonempty 
set $X$.  
Then the structure group $\CG$ is $2$-step nilpotent.
\end{proposition}

\begin{proof}
The group $G\type 3$ acts on $X\type 3$ in the natural way. Let 
$g\in\CG$ and $\ubx\in\CQ$. Write $\ubx=(\bx',\bx'')$ with $\bx, \bx'\in 
\CP$. 
Then $g\type 2\cdot\bx'$ is a parallelogram equivalent to $\bx'$.  
Thus $(g\type 2\cdot\bx',\bx')\in\CQ$ and so $(g\type 
2\cdot\bx',\bx'')\in\CQ$. By symmetries of $\CQ$,
for every face $f$ of $\{0,1\}^3$ we have that $g\type 
{3,f}\cdot\bx\in\CQ$. Therefore, $\CQ$ is invariant under the 
group $G\type{3,2}$.

Let $g\in \CG_3$. By Lemma~\ref{lem:comut}, part~\ref{it:comut5}
$(1,\dots,1, g)\in G\type {3,2}$.  For every $x\in X$, we have
 $(x,\dots,x,x)\in\CQ$ and 
thus $(x,\dots,x,g\cdot x)\in\CQ$. 
Since $\CQ$ is a strong 
structure, $g\cdot x=x$ and $g=1$.  Thus $\CG$ is $2$-step nilpotent.  
\end{proof}

In all examples considered thus far, 
the group $G$ is included in $\CG(X)$. 
But we note that $\CG(X)$ may be substantially larger than $G$.
Consider the situation of Example~\ref{ex:abel2}. Let 
$\phi\colon B\to F$ be 
a group homomorphism and define a transformation $h$ of $X$ by
$$
h\cdot x= \phi(\pi(x))\cdot x\ .
$$
Then $h$ belongs to $\CG(X)$ and is not translation by an element 
of $G$.

\section{Description of parallelepipeds structures}
\label{sec:desc}
Henceforth, $(\CP,\CQ)$ is a (strong) parallelepiped structure on a set $X$.

Our goal is to characterize when a parallelepiped structure 
is a nilparallelepiped structure, and we do so in Theorem~\ref{th:structure}.  
Moreover,  in Theorem~\ref{th:split} and Corollary~\ref{th:embed}, 
we give sufficient conditions for this property to hold.  We also 
give an example of a parallelepiped structure without this property 
(Example~\ref{ex:counterex}) and show that in the general case a 
parallelepiped structure can be imbedded in a nilstructure 
(Proposition~\ref{prop:embedd}) in a sense 
explained below.

For parallelograms we use the notation introduced in 
Section~\ref{sec:moreexamples}.  
$B$ is the base group of the parallelogram 
structure $\CP$ and $\pi\colon X \to B$ 
the surjection defined in Section~\ref{subsec:def-paragrammes}.
 We recall that for $x,y\in X$, $\wec{x,y}$ is the equivalence class 
 of the pair $(x,y)$ under the relation $\simp$, that is
 $$
 \wec{x,y}=\pi(y)\,\pi(x)\inv\ .
 $$
For $b\in B$, the \emph{fiber} $F_b$ of $b$ is defined by 
$F_b:=\pi\inv(\{b\})$. 
 
For parallelepipeds we use the notation of 
Section~\ref{subsec:def-parapedes}.
The equivalence class (under the relation $\approx$) of a 
parallelogram $\bx$ is denoted by $\br{\bx}$.   
We denote the quotient space $\CP/\approx$ by $P$.

\subsection{The groups $P_s$}
\label{sec:Ps}
For every $s\in B$, we define
\begin{align*}
    X_s & :=\bigl\{(x_{01}, x_{10})\in X^2\colon \wec{x_{01},x_{02}}=s\}
\subset X^2 \ .\\
\CP_s  & := \bigl\{ \bx\in\CP\colon 
\wec{x_{00},x_{01}}=s\bigr\}\subset X\type 2\ .\\
\CQ_s & := \bigl\{ \uby\in\CQ\colon 
\wec{y_{000},y_{001}}=s\bigr\}\subset X\type 3\ .
\end{align*}
If two parallelograms $\bx$ and $\by\in\CP$ are equivalent under
the relation $\approx$, then
 $\wec{x_{00},x_{01}}= \wec{y_{00},y_{01}}$ and thus they belong to the 
same $\CP_s$. 
Therefore each set $\CP_s$ is a union of equivalence classes under 
the relation $\approx$.  
Writing $P_s$ for the set of equivalence classes of parallelepipeds belonging 
to $\CP_s$, we have a partition of $P$:
$$ 
P=\bigcup_{s\in B}P_s\ .
$$

We identify $X\type 2$ with $(X\times X)^2$ in the natural way: 
$$
\bigl(x_{00},x_{01},x_{10},x_{11}\bigr) =   
\bigl((x_{00},x_{01}),(x_{10},x_{11})\bigr)
$$
and $X\type 3$ with $(X^2)\type 2$:
\begin{multline*}
     \bigl(
x_{000},x_{001},x_{010},x_{011},x_{100},x_{101},x_{110},
x_{111}
\bigr)\hfill\qquad\qquad\strut
\\
= 
\bigl(
(x_{000},x_{001}),(x_{010},x_{011}),(x_{100},x_{101}),(x_{110},
x_{111})
\bigr)\ .
\end{multline*}
Therefore we  view $X_s\times X_s$ as a subset of $X\type 2$ and 
$X_s\type 2$ as a subset of $X\type 3$.  
We have:
$$
\CP_s =\CP\cap(X_s\times X_s)\text{ and }
\CQ_s= \CQ\cap X_s\type 2\ .
$$
We reformulate this for clarity. A pair of elements of $X_s$ 
represents four points in $X$ forming 
a parallelogram of $X$ and this parallelogram belongs to $\CP_s$.
An element  $\ubx\in X_s\type 2$ consists in four points of $X_s$, 
that is, in eight points of $X$, and these eight points form a 
parallelepiped in $X$ if and only if $\ubx\in\CQ_s$. 

\begin{lemma}
\label{lemma:Qs}
Maintaining the above notation, 
$\CQ_s$ is a parallelogram structure on $X_s$.
\end{lemma}

\begin{proof}
All the needed properties are immediate other than 
transitivity.
Let $(x_0,x_1)$, $(x_2,x_3)$, $(x_4,x_5)$ be three paints in $X_s$. 
Choose $x_6\in X$ with $\wec{x_4,x_6}=\wec{x_0,x_2}$. 
The seven points $x_0,\dots,x_6$ satisfy the closing parallelepiped 
property and so
there exists $x_7\in X$ such that $(x_0,\dots,x_6,x_7)\in\CQ$. 
This parallelepiped actually belongs to $\CQ_s$ and thus also
$\bigl((x_0,x_1),(x_2,x_3),(x_4,x_5),(x_6,x_7)\bigr)$ and our claim 
is proved.
\end{proof}

We note that this parallelogram structure is not strong, as there is 
freedom in the choice of $x_6$ in the above construction.

Recall that $\simqs$ denotes the equivalence relation on $X_s^2$ 
associated to the parallelogram structure $\CQ_s$ on $X_s$:
two pairs  of points in $X_s$ are equivalent for this relation 
if the form a parallelogram in $\CQ_s$. 
If we consider these two pairs as parallelograms in $X$, 
then these parallelograms belong to $\CP_s$ and 
Lemma~\ref{lemma:Qs} implies that they are equivalent 
under the relation $\approx$. 
Therefore we can identify the two quotient spaces
$$
(X_s)^2/\!\simqs\quad =\quad \CP_s/\!\approx\quad =\quad  P_s\ .
$$
By Lemma~\ref{lem:mult-B}, the quotient space $(X_s)^2/\simqs=P_s$  
can be endowed with a multiplication, that gives it the structure of an abelian 
group.
For clarity, this multiplication in 
the present notation, viewing $P_s$ as $\CP_s/\approx$.

Let $u,v$ be two classes in $P_s$.
Let $(x_0,x_1,x_2,x_3)$ be a parallelogram in the class $u$.
As $\wec{x_2,x_3}=s$ there exist two points $x_4$ and $x_5$ in $X$ 
such that $(x_2,x_3,x_4,x_5)$ is a parallelogram in the class $v$. 
Then $uv$ is the class of the parallelogram
$(x_0,x_1,x_4,x_5)$.
\subsection{A homomorphism}
\label{sec:qs}
Let $s\in B$.  
If two parallelograms $\bx$ and $\by$ of $X$ belonging to $\CP_s$ are equivalent 
under the relation $\approx$, then
 $\wec{x_{00},x_{10}}=\wec{y_{00},y_{10}}$.
Therefore there exists a map $q_s\colon P_s\to B$ such that
$$
q_s(\br{\bx})= \wec{x_{00},x_{10}}
\text{  for every parallelogram }\bx\in \CP_s\ .
$$
This map is clearly a group homomorphism from $P_s$ onto $B$.
The kernel of this homomorphism consists in the set of equivalence 
classes of parallelograms $\bx\in\CP$ with $\wec{x_{00},x_{01}}=s$ and 
$\pi(x_{10})=\pi(x_{00})$.

\section{The fiber group}
In this Section, $(\CP,\CQ)$ is a strong parallelepiped structure on 
a nonempty set $X$ and we maintain the notation of the preceding 
Section.

\subsection{Vertical parallelograms}
We begin with some simple observations and some more 
vocabulary. 
\begin{lemma}
\label{lemma:double} \strut
\begin{enumerate}
\item
\label{cl:double} 
If $\bx\in\CP$, then 
$\br{x_{00},x_{00},x_{01},x_{01}}=
\br{x_{10},x_{10},x_{11},x_{11}}$ and \\
$\br{x_{00},x_{01},x_{00},x_{01}}=
\br{x_{10},x_{11},x_{10},x_{11}}$.

\item
\label{cl:samefiber}
If $x$ and $y$ belong to the same fiber, then 
$\br{x,x,y,y}=\br{x,y,x,y}=1$.
\end{enumerate}
\end{lemma}

\begin{proof}
For part~\ref{cl:double}, since $(\bx,\bx)\in\CQ$, both properties follow  from the symmetries of 
$\CQ$.
By assumption, $(x,x,x,y)\in\CP$ and part~\ref{cl:samefiber}
follows from part~\ref{cl:double}.
\end{proof}

If $x_{00},x_{01}, x_{10},x_{11}$ are four points in the same fiber, then 
$(x_{00},x_{01}, x_{10},x_{11})\in\CP$.  Thus it makes sense to define:
\begin{definition}
A parallelogram with $4$ points in the same fiber is 
called a \emph{vertical parallelogram}.  
\end{definition}

A parallelogram equivalent to a vertical one is also vertical, and 
thus the family of vertical parallelograms is a union of equivalence classes. 

\begin{lemma}
\label{cl:croise}
If 
$x_{00},x_{01}, x_{10},x_{11}$ are four points in the same fiber, 
then
$$\br{x_{00},x_{01},x_{10},x_{11}}=\br{x_{00},x_{10},x_{01},x_{11}}\ .$$
\end{lemma}
\begin{proof}
There exists a unique $y\in X$ such that 
$\br{x_{00},x_{01}, x_{10},x_{11}} = 
\br{x_{00},x_{01},x_{01},y}$.  Thus 
$(x_{00},x_{01}, x_{10},x_{11},x_{00},x_{01},x_{01},y)\in\CQ$. 
By symmetry, 
$(x_{00},x_{10},x_{01},x_{11},x_{00},$ $x_{01},x_{01},y)\in\CQ$  and 
thus $\br{x_{00},x_{10},x_{01},x_{11}}=
\br{x_{00},x_{01},x_{01},y}=\br{x_{00},x_{01},x_{10},x_{11}}$.
\end{proof}

\subsection{The fiber group $F$ and its action on $X$}
\begin{notation}
We denote the unit element of $B$ by $1$ and  let $F$ denote the 
kernel of the homomorphism $q_1\colon P_1\to B$.
\end{notation}

In other words, $F$ is the set of equivalence classes (under the relation 
$\approx$) of vertical parallelograms. Recall that the 
multiplication in $F$ satisfies:
\begin{multline}
\label{eq:multF}
\text{if } x_0,x_1,x_2,x_3,x_4,x_5\text{ belong to the same fiber, then}\\
\br{x_0,x_1,x_2,x_3}\, \br{x_2,x_3,x_4,x_5}= \br{x_0,x_1,x_4,x_5}\ .
\end{multline}

We now define an action of $F$ on $X$, mapping each fiber to itself and 
use this to describe the vertical parallelograms. 

Let  $x\in X$ and $u\in F$. Recall that $u$ is the class 
of some vertical parallelogram. It thus follows from the closing 
parallelepiped property that there exists $y\in X$ such that
$\br{x,x,x,y}=u$.   As the structure is strong, 
this point $y$ is unique.

\begin{notation}
For $x\in X$ and $u\in F$, we write $u\cdot x$ for the point of $X$ 
defined by $\br{x,x,x,y}=u$.
\end{notation}
\begin{lemma}
\label{lem:actionF}
The map $(u,x)\mapsto u\cdot x$ is an action of the group $F$ on the 
set $X$.
This action preserves each fiber and acts transitively and freely on 
each fiber.
\end{lemma}
\begin{proof}
Note that $1_P\cdot x=x$ for every $x\in X$. 
We are left with showing that for $u,v\in F$ and $x\in X$, we have 
$(vu)\cdot x=v\cdot(u\cdot x)$. 
Let $y=u\cdot x$ and $z=v\cdot y=v\cdot(u\cdot x)$. 
Then $v = \br{u\cdot x, u\cdot x, u\cdot x, v\cdot(u\cdot x)} = 
\br{x, u\cdot x, x, u\cdot x}\br{x, u\cdot x, x, v\cdot(u\cdot x)} = 
\br{x,u\cdot x,x,v\cdot(u\cdot x)}$ and 
$\br{x,x,x,v\cdot(u\cdot x)}
=\br{x,x,x,u\cdot x}\,\br{x,u\cdot x,x,v\cdot(u\cdot x)}=uv=vu$
because $F$ is abelian.  Thus we have an action.  
 
By construction, $F$ preserves each fiber. For every $x,y$ in the 
same fiber, there exists a unique $u\in F$ such that $u\cdot x=y$, 
namely $u=\br{x,x,x,y}$. This means that the action of $F$ on each 
fiber is free and transitive.
\end{proof}

\begin{proposition}
\label{prop:vertical}\strut
\begin{enumerate}
\item
\label{it:eqclass}
Let $(x,u\cdot x,v\cdot x, w\cdot x)$ be a vertical parallelogram. Then the 
equivalence class in $F$ of this parallelogram is the element $wu\inv v\inv$ of $F$.

\item 
\label{it:pareq}
For every parallelogram $\bx\in\CP$ and $u\in F$, we have $(u\cdot x_{00},u\cdot x_{01}, x_{10},x_{11})\in\CP$
and $\br{u\cdot x_{00},u\cdot x_{01}, x_{10},x_{11}}=\br{\bx}$.

\item 
\label{cl:F21}
Every transformation in $F\type{2,1}$ maps every parallelogram to an equivalent one.
\end{enumerate}
\end{proposition}

\begin{proof}
We have 
\begin{align*}
w = &\br{x,x,x,w\cdot x}=
\br{x,x,x,u\cdot x}\,
\br{x,u\cdot x,v\cdot x,w\cdot x}\,\br{v\cdot x,w\cdot x,x,w\cdot x}\\
 = &u\,\br{x,u\cdot x,v\cdot x,w\cdot x}\,\br{v\cdot x,x,w\cdot x,w\cdot x}\\
 = &u\,\br{x,u\cdot x,v\cdot x,w\cdot x}\,\br{v\cdot x,x,x,x}\,\br{x,x,w\cdot x,w\cdot x} \\
 = &u\,\br{x,u\cdot x,v\cdot x,w\cdot x}\,\br{v\cdot x,x,x,x}
=u\,\br{x,u\cdot x,v\cdot x,w\cdot x}\,\br{x,x,x,v\cdot x}\\
 = &u\,\br{x,u\cdot x,v\cdot x,w\cdot x}\,v\ .
\end{align*}
This proves part~\ref{it:eqclass}.  

We now prove part~\ref{it:pareq}.  Let $\bx\in\CP$ 
and let $u\in F$. By Lemma~\ref{lemma:double}, part~\ref{cl:double} 
we have 
$\br{x_{00},x_{01},x_{00},x_{01}}=\br{x_{10},x_{11},x_{10},x_{11}}$.  
By part~\ref{it:eqclass}, 
$\br{x_{00},x_{00},u\cdot x_{00},x_{00}}
=u\inv=\br{x_{01},x_{01},u\cdot x_{01},x_{01}}$
and thus
$\br{x_{00},x_{01}, u\cdot x_{00},u\cdot x_{01}}=\br{x_{00},x_{01},x_{00},x_{01}}=
\br{x_{10},x_{11},x_{10},x_{11}}$ by Lemma~\ref{lemma:double}, part~\ref{cl:double}.  
The claim follows by symmetry.

The group $F\type 2$ acts on $X\type 2$ coordinate-wise.  
By part~\ref{it:pareq} and the symmetries of $\CQ$, we have 
the statement in part~\ref{cl:F21}.  
\end{proof}

>From this lemma, we deduce that 
$\CQ$ is invariant under the subgroup $F\type{3,1}$ of $F\type 3$.
Conversely, we have:
\begin{proposition}
\label{cl:F21bis}
Let $\bx\in\CP$ and $\bu\in F\type 2$. 
If the parallelograms $\bx$ and $\bu\cdot\bx$ are equivalent,   
then $\bu\in F\type{2,1}$.
\end{proposition}

\begin{proof}
There exists $u'_{11}$ such that $(u_{00},u_{01},u_{10},u'_{11})\in F\type{2,1}$ and 
$\br{u_{00}\cdot x_{00},u_{01}\cdot x_{01},u_{10}\cdot x_{10},u'_{11}\cdot x_{11}}=\br{\bx}=
\br{u_{00}\cdot x_{00},u_{01}\cdot x_{01},u_{10}\cdot x_{10},u_{11}\cdot x_{11}}$. 
Since the structure is strong, $u'_{11}\cdot x_{11}=u_{11}\cdot x_{11}$. 
But since $F$ acts freely on each fiber, $u'_{11}=u_{11}$.
\end{proof}

\begin{remark}
Returning to the parallelogram structure $\CQ_s$ on $X_s$, we have 
already noticed that it is not strong.  
As in Section~\ref{subsec:def-paragrammes}, we have a strong  
parallelogram structure by taking a quotient of $X_s$ by some 
equivalence relation.  Using the preceding proposition, we have that 
two elements $(x,y)$ and $(x',y')$ of $X_s$ are equivalent for this 
relation if and only if there exists $u\in F$ with $x'=u\cdot x$ and 
$y'=u\cdot y$.
\end{remark}

\begin{proposition} 
The group $F$ is included in the center of $\CG$.
\end{proposition}

\begin{proof}
For every $u\in F$, we have $u\type 2\in F\type{2,1}$ and so $F\subset\CG$.

Let $g\in\CG$ and $u\in F$. 
Let $f$ be a face of the cube and $e$ an edge of the cube with $f\cap e=\{111\}$. 
Then the transformations $g\type{3,f}$ and $u\type{3,e}$ map 
$\CQ$ to itself, and thus so does the commutator of these transformations. 
It is immediate to check that
this commutator is equal to $[g;u]\type{3,e}$.   
This means that for $\bx\in\CQ$, we also have 
$(x_{000},x_{001},\dots,x_{110}, [g;u]\cdot x_{111})\in\CQ$. 
As the structure is strong, 
$[g;u]\cdot x_{111}=_{111}$. 
Thus $[g;u]$ is the identity transformation.
\end{proof}

\subsection{Structure Theorem}
We are now ready to characterize parallelepiped structures 
that are nilstructures.  Recall (see Section~\ref{subsec:def-paragrammes}) 
that  $i$ is point in 
$X$ and that the map $\pi\colon X\to B$
is defined by $\pi(x)=\ov{i,x}$ for every $x\in X$.

\begin{theorem}
\label{th:structure}
Let $G$ be a subgroup of $\CG$ containing $F$ and 
assume that $G$ acts transitively
on $X$. Let $\Gamma$ be the stabilizer of
some $i\in X$ in $G$ and 
identify $X$ with $G/\Gamma$ in the natural way.
Then the parallelogram structure $(\CP,\CQ)$ on $X$ coincides with
the structure $(\CP_X,\CQ_X)$ associated to $G,F$ and $\Gamma$ as in
Proposition~\ref{prop:nil}
\end{theorem}

Thus the parallelepiped structure $(\CP,\CQ)$  is isomorphic 
to a nilmanifold parallelepiped 
structure as in Proposition~\ref{prop:nil} 
if and only if its structure group $\CG$ acts transitively on $X$.

\begin{proof}
\emph{Step 1.}
We first show that $G,F$ and $\Gamma$ satisfy the hypotheses of 
Proposition~\ref{prop:nil}.

We have already shown that $F$ is included in the center of $G$. Since
$F$ acts freely on $X$, $\Gamma\cap F=\{1\}$.  We are left with showing:
\begin{claim}
$F$ contains  the commutator subgroup $\CG_2$ of $\CG$.
\end{claim}

Recall that $\CG\type{3,2}\supset\CG_2\type{3,1}$. 
Therefore every element of $\CG_2\type{3,1}$ leave $\CQ$ invariant and thus 
every element of $\CG_2\type{2,1}$ maps each parallelogram to an equivalent
one.
On the other hand, since $B$ is abelian $g$ belongs to the kernel of $p$ 
and thus maps every point of $X$ to a point in the same fiber.

Let $g\in\CG_2$ and $x,y\in X$. There exist $u,v\in F$ with $g\cdot x=u\cdot x$ 
and $g\cdot y=v\cdot y$. As $(x,y,x,y)\in\CP$ and $(g,g,1,\tau)\in \CG_2\type{2,1}$ we have 
$\br{x,y,x,y}=\br{g\cdot x,g\cdot y,x,y}=\br{u\cdot x,v\cdot y,x,y}$. 
By Proposition~\ref{cl:F21bis}, we have $u=v$. 

Therefore there exists $u\in F$ with $g\cdot x=u\cdot x$ for every $x$, meaning that  $g=u$.
This proves the claim.  
Therefore we can define the structure $(\CP_X,\CQ_X)$ as in 
Proposition~\ref{prop:nil}.

\noindent\emph{Step 2.}
Let $g\in\CG$. For arbitrary $x,y\in X$, $g\type 2$ maps the parallelogram $(x,x,y,y,)$ to an equivalent one 
and this implies that 
$\pi(g\cdot x)\,\pi(x)\inv=\pi(g\cdot y)\,\pi(y)\inv$.
 Therefore there exists an element $p(g)\in B$ such that
$$
\pi(g\cdot x)= p(g)\,\pi(x)\text{ for  every }  x\in X\ .
$$
The map $p\colon G\to B$ defined in this way is clearly a group
 homomorphism.
We show that 
\begin{claim}
$\ker(p)=F\Gamma$.
\end{claim}

The kernel of $p$ contains clearly contains  $F$.
For $\gamma\in\Gamma$ we have $\pi(i)=\pi(\gamma\cdot i)
=p(\gamma)\,\pi(i)$ and thus $\gamma\in\ker(p)$. 
We get that $\ker(p)\supset F\Gamma$.
Conversely, let $g\in\ker(p)$. Then $g\cdot i$ belongs to the same 
fiber as $i$ 
and there exists $u\in F$ with $g\cdot i=u\cdot i$ and we have 
$u\inv g\in\Gamma$.
This proves the claim.  

As the parallelogram structure $\CP$ is associated to the projection $X\to B$ and 
$\CP_X$ is associated to the projection $X\to \CG/F\Gamma$ we have that 
$$
\CP=\CP_X\ .
$$
\noindent\emph{Step 3.} We are left with showing that $\CQ=\CQ_X$.

Since$(i,i,\cdots,i)\in\CQ$, for every $\bg\in G\type{3,2}$ we have 
$(g_{000}\cdot i,g_{001}\cdot i,\cdots,g_{111}\cdot i)\in\CQ$ by definition of the 
group $\CG$. Under our identification it means that $\CQ_X\subset\CQ$.

Let $\ubx\in\CQ$. As $\CP=\CP_X$, each face of $\ubx$ belongs to 
$\CP_X$ and by  condition~\ref{rem:closing2} (applied to the structure 
$(\CP_X, \CQ_X)$)  
there exists $x'_{111}$ such that 
$(x_{000},x_{001},\dots,x_{110},x'_{111})\in \CQ_X$. As $\CQ_X\subset\CQ$ these 
$8$ points form a parallelepiped in $\CQ$. As this last structure is strong,
 $x'_{111}=x_{111}$ and thus $\ubx\in\CQ_X$.
\end{proof}

\section{Conditions for transitivity}
\label{sec:trans}

For  $s\in B$, let $F_s$ be the kernel of the group homomorphism 
$q_s\colon P_s\to B$ defined in Section~\ref{sec:qs}.  
Thus $F_s$ is the family 
of $\approx$-equivalence classes of parallelograms $\bx$ with 
$\wec{x_{00},x_{01}}=s$ and $\wec{x_{00},x_{10}}=1$.  This means that
$x_{00}$ and $x_{10}$ lie in the same fiber and therefore $x_{01}$ and 
$x_{11}$ also lie the same fiber.

Therefore, if $\bx$ is a parallelogram its equivalence class $\br{\bx}$ belongs 
to $F_s$ if and only if $\bx$ can be written under the form:
$$ 
(a,b, u\cdot a,v\cdot b)\text{ with } a,b\in X,\ \wec{a,b}=s
\text{ and } u,v\in F\ .
$$
Let $a',b'\in X$ with $\wec{a',b'}=s$ and $u',v'\in F$. 
Then 
$\br{a,b, u\cdot a,v\cdot b}=\br{a',b', u'\cdot a',v'\cdot b'}$
if and only if
$\br{a,b,a',b'}=\br{u\cdot a,v\cdot b, u'\cdot a',v' \cdot b'}$.
By Proposition~\ref{prop:vertical}, part~\ref{cl:F21} this last class is equal to 
$\br{a,b,a', uv\inv {u'}\inv v'\cdot b'}$.  This is equivalent to
$\br{a,b,a',b'}$ if and only if $uv\inv {u'}\inv v'=1$.

Therefore by passing to the quotient, the map
$$
(a,b, u\cdot a,v\cdot b)\mapsto vu\inv
$$
induces a one to one map $j_s\colon F_s\to F$.

The map $j_s$ is clearly onto. Moreover, it follows immediately from the 
definition of the multiplication in $F$ that this map is a group 
homomorphism. We conclude that $j_s\colon F_s\to F$ is a group 
isomorphism.

Identifying $F$ with $F_s$ through this isomorphism, for every $s\in B$ 
we have an exact sequence:
\begin{equation}
    \label{eq:exact}
0\longrightarrow F\longrightarrow P_s\longrightarrow B
\longrightarrow 0\ .
\end{equation}

\begin{theorem}
\label{th:split}
The group $\CG$ acts transitively on $X$ if and only if  for every $s\in B$, the exact 
sequence~\eqref{eq:exact} splits.  
\end{theorem}

\begin{proof}
Since the subgroup $F$ of $\CG$ acts transitively on each fiber, 
then the group $\CG$ acts transitively on $X$ if and only if the map 
$q\colon \CG\to B$ is onto.

\medskip\noindent\emph{First step.} Let $s\in B$.
First we assume that the exact sequence~\eqref{eq:exact} splits.
This means that there 
exists a group homomorphism $\phi\colon B\to P_s$ with $q_s\circ \phi=\id_B$.
We build a transformation $h$ of $X$, belonging to $\CG$, with $q(g)=s$.

 We choose $a,b\in X$ with $\wec{a,b}=s$.
 Let $x\in X$. Then we define 
$g\cdot x$ to be
the unique point $g\cdot x$ in $X$ such that $(a,b,x,g\cdot x)$ is a 
parallelogram in the class $\phi(\wec{a,x})$.

Let $\bx$ be a parallelogram in $X$. We have
$\wec{a,x_{00}}\,{\wec{ a,x_{01}}}\inv\,
{\vec {ax_{10}}}\inv\,\wec{a,x_{11}}=1$. 
Since $\phi$ is a group homomorphism, we have 
$$
\phi(\wec{a,x_{00}})\,\phi({\wec{ a,x_{01}}}\inv)\,
\phi({\vec {ax_{10}}}\inv)\,\phi(\wec{a,x_{11}})=1\ ,$$ 
This means that 
$$
\br{ a,b,x_{00},g\cdot x_{00}}\,
\br{ a,b,x_{01},g\cdot x_{01}}\inv\,
\br{ a,b,x_{10},g\cdot x_{10}}\inv\,
\br{ a,b,x_{11},g\cdot x_{11}}=1 \ .
$$
By the definition of the multiplication in $P_s$, we have that
$\br{x_{00},g\cdot x_{00}, x_{01}, ,g\cdot x_{01}}=
\br{x_{10},g\cdot x_{10}, x_{11}, ,g\cdot x_{11}}$
and thus $\br{\bx}=\br{g\cdot x_{00},g\cdot x_{01},
g\cdot x_{10}, g\cdot x_{11}}$.

Therefore the transformation $g$ of $X$ maps every parallelogram to an 
equivalent one and thus it belongs to $\CG$. 
By construction, $q(g)=s$.

\medskip\noindent\emph{Second step.}
Conversely, let $s\in B$ and assume that there exists $g\in \CG$ with $q(g)=s$.

We choose $a\in X$ and define $b=g\cdot a$.
For every $x\in X$, define $\psi(x)\in P_s$ to be the 
$\approx$-equivalence class of the parallelogram $(a,b,x,g\cdot x)$.
Then for every $x$, we have  $q_s(\psi(x))=\wec{a,x}$. Moreover, 
$\psi(a)=\br{a,b,a,b}=1$ and for every $x\in X$ we have 
$q_s(\psi(x))=\wec{a,x}$.

By the definition of $\CG$ and 
the same computation as in the first step, we have that 
\begin{equation}
\label{eq:parapsi}
\text{when $\bx$ is a 
parallelogram, } 
\psi(x_{00})\psi(x_{01})\inv\psi(x_{10})\inv\psi(x_{11})=1\ .
\end{equation}

We deduce that for $x\in X$,  $\psi(x)$ depends only in 
$\pi(x)$ and thus on $\wec{a,x}$.  Thus
there exists a map $\phi\colon B\to P_s$ with
$\psi(x)=\phi(\wec{a,x})$ for every $x\in X$.

 The  relation~\eqref{eq:parapsi}  thus implies 
that, when $b_{00}, b_{01}, b_{10}, b_{11}$  are four points in $B$ 
with $b_{00} b_{01}\inv b_{10}\inv b_{11}=1$ then 
$\phi(b_{00})\phi(b_{01})\inv\phi(x_{10})\inv\phi(x_{11})=1$.
As $\phi(1)=\psi(a)=1$, it follows that $\phi$ is a group homomorphism.

We have $q_s\circ \phi=\id_B$ and so the exact 
sequence~\eqref{eq:exact} splits.
\end{proof}

It follows immediately that: 
\begin{corollary}  
\label{th:embed}
If either the fiber group is injective or if the 
base group is projective, then the parallelepiped structure 
is a nilparallelepiped structure, as in Proposition~\ref{prop:nil}.
\end{corollary}

\subsection{Imbeddings}

In this Section we first show that there exist parallelepiped structures 
such that the structure group does not act transitively 
on the space, and so the parallelepiped structure is not a 
nilstructure.  After the example, 
we show that every parallelepiped structure can be imbedded in a 
nilstructure.  We start with a proposition that defines this imbedding.  
\begin{proposition}
\label{prop:embedd}
Let $(\CP,\CQ)$ be a parallelepiped structure on a set $X$, assume 
that $B$ is the base is the base of $X$, $\pi\colon X\to B$ 
the natural projection, and let $Y$ be 
a subset of $X$ satisfying the following two conditions:
\begin{enumerate}
\item
The restriction to $Y$ of $\pi$ is onto.
\item
If $\ubx\in\CQ$ is a parallelepiped in $X$ such that the seven points
$x_{000},\dots,x_{110}$ belong to $Y$, then $x_{111}\in Y$.
\end{enumerate}
Then $(\CP\cap Y\type 2,\CQ\cap Y\type 3)$ is a parallelepiped 
structure on $Y$.
\end{proposition}

\begin{definition}
If the conditions in the proposition 
are satisfied, we say that the parallelepiped structure 
$(Y,\CP\cap Y\type 2,\CQ\cap Y\type 3)$ is \emph{embedded} in 
$(X,\CP,\CQ)$.
\end{definition}

The proof of the proposition follows immediately from the definitions.  

\begin{example}[A parallelepiped structure which is not a nilstructure]
\label{ex:counterex}
Let $ p>2$ be a prime number and $B=\Z/p\Z$ and $F=\Z/p^2\Z$.

In this example we use additive notation for $B$ and $F$.
Let $(\CP,\CQ)$ the (abelian) parallelepiped structure on 
$X=B\times F$ defined as in 
Example~\ref{ex:abel2}, but here we prefer additive notation:
$$
\CP= B\type{2,1}\times F\type 2\ , \ 
\CQ=B\type{3,2}\times F\type{3,1}\ .
$$

The base group is $B$ and $\pi\colon X\to B$ is projection onto 
the first coordinate.

In this case the groups $P_s$, $s\in B$, introduced in Section~\ref{sec:Ps} 
can easily be defined explicitly and we do so now. 

Let $s\in B$. Then $P_s=B\times F$. The element  
$\bigl((b_{00}, f_{00}),((b_{01}, f_{01}),((b_{10}, f_{10}),((b_{11}, 
f_{11})\bigr)$ of $X\type 2=B\type 2\times F\type 2$ 
is a parallelogram belonging to $\CP_s$ if and only if
$b_{01}-b_{00}=b_{11}-b_{10}=s$.  In this case, the class of this 
parallelogram is the element $(b_{10}-b_{00}, 
f_{00}-f_{01}-f_{10}+f_{11})$ of $P_s$.
The homomorphism $q_s\colon P_s\to B$ is the first coordinate and 
the imbedding $F\to P_s$ is the map $f\mapsto (0,f)$.

Let $X'$ be the subset of $X$ given by
$$
X'=\bigl\{ (b,f)\colon f=\frac {b(b-1)}2\bmod p\Z/p^2\Z\bigr\}\ .
$$
It is easy to check that if $\bx$ is a parallelepiped of $X$ such that 
seven of its edges belong to $X'$, then the remaining edge also belongs to 
$X'$. Therefore, $\CP':=\CP\cap {X'}\type 2$ and 
$\CQ':=\CQ\cap {X'}\type 3$ is a parallelepiped structure on $X'$, 
imbedded in the parallelepiped structure on $X$.  

The fiber group is $F'=p\Z/p^2\Z\cong\Z/p\Z$. 
For $s\in B$, the family of parallelograms of this structure is 
written $\CP'_s$ that is, $\CP'_s=\CP_s\cap \CP'$.
Let $P'_s$ be the abelian group of classes of 
parallelograms of this family. An immediate computation shows that 
$P'_s$ is the subgroup
$$
P'_s=\bigl\{ (b,f)\in B\times F\colon f=sb\bmod p\Z/p^2\Z\bigr\}
$$
of $P_s$. For $\neq 0$ this group is isomorphic to $\Z/p^2\Z$ and 
thus it is not isomorphic to the direct sum $B'\oplus F'\approx
\Z/p\Z\oplus\Z/p\Z$. Therefore the exact sequence~\eqref{eq:exact} 
does not split and so there is no $g\in\CG(X)$ projecting on 
$s$.
\end{example}

\begin{theorem}
\label{th:nik-imbedd}
Every parallelepiped structure can be imbedded in a nilparallelepiped 
structure.  
\end{theorem}

\begin{proof}
Let $(\CP,\CQ)$ be a parallelepiped structure on the set $X$, let 
$B$ be the base group and let $F$ be the fiber group.  

The abelian group $F$ can be imbedded as a subgroup of a divisible group $E$.
We write $X\otimes E$ for the set $X\times E$, quotiented by the 
equivalence relation given by
$$
(x,e)\cong (x',e')\text{ if there exist } u\in F\text{ with }
x'=u\cdot x\text{ and } e'=u\inv e\ .
$$
We write $j\colon X\times E\to X\otimes E$ for the quotient map.

We now define a parallelepiped structure $(\CP_E,\CQ_E)$ on $X\otimes E$.

Let $(x,e), (x',e')\in X\times E$. 
If $(x,e)\cong (x',e')$, then $\pi(x)=\pi(x')$. 
Therefore we can define a map $\pi_E\colon X\otimes E\to 
B$ by mapping the equivalence class  $j(x,e)$ to $\pi(x)$.
We define $\CP_E$ to be the parallelogram structure on $X\otimes E$ 
associated to this projection.

Let $(\CP_E,\CQ_E)$ be the parallelepiped structure on the abelian group 
$E$ defined as in 
Example~\ref{ex:abel1} 
and let  $\CQ_{X\otimes E}$ be the image of 
$\CQ\times\CQ_E$ under  
the natural projection $X\type 3\times E\type 3= (X\times E)\type 3\to (X\otimes E)\type 3$.
We claim that $(\CP_E,\CQ_E)$ is a parallelepiped structure
on $X\otimes E$.

For every $\bx\in\CQ_{X\otimes E}$, every face of $\bx$ obviously 
belongs to $\CP_{X\otimes E}$. 
The symmetries of $\CQ_{X\otimes E}$ are 
obvious. 
Transitivity  follows from Proposition~\ref{prop:vertical}, part~\ref{cl:F21} 
and Proposition~\ref{cl:F21bis}.   
The closing parallelepiped property follows in the same way that 
Lemma~\ref{cl:closing3} does.  

The base group of $\CP_{X\otimes E}$ is $B$ and 
the fiber group of  $\CQ_{X\otimes E}$ is $F$. Since $F$ is a divisible 
group, the structure $(\CP_{X\otimes E},\CQ_{X\otimes E})$ on 
$X\otimes 
E$ is a nilstructure.

On the other hand, the map $X\to X\otimes E$ associating $x\in X$ to the
equivalence class of $j(x,1)$ is one to one,  and we can consider $X$ as a 
subset of $X\otimes E$.
Thus $\CP=\CP_{X\otimes E}\cap X\type 2$ and 
$\CQ=\CQ_{X\otimes E}\cap X\type 3$.
\end{proof}

\section{Higher levels?}
\label{sec:higher}

Gowers norms (for all $k\geq 2$) have already 
been used in several contexts and it is natural to ask to 
what extent the results here can be generalized for $k\geq 4$.

In the setting of ergodic theory, the authors define 
seminorms for all $k\geq 1$ and measures on Cartesian products of the 
space playing the same role played by the structures of this paper.
The ``strong'' structures of this paper correspond to the ``systems 
of order $k$'' of~\cite{HK}. These systems are completely  
characterized in terms of nilmanifolds.  But the descriptions we 
give in the context of the present paper are substantially weaker.

Our definitions of parallelogram and parallelepiped structures 
extend immediately to structures of any dimension $k$, for which 
basic models are given by $(k-1)$-step nilmanifolds. 

The main constructions and results of Section~\ref{sec:parapedes} extend 
directly, for example the reduction 
to strong structures (Proposition~\ref{prop:weakstrong}),
the definition of the seminorm (Proposition~\ref{prop:norme3}), and the 
definition of the structure group (Definition~\ref{def:struc}).  
In particular, this structure group is always $(k-1)$-step nilpotent.

The results of Sections~\ref{sec:desc}-\ref{sec:trans} 
are more difficult to extend.  The main difficulty is that 
a good description of a structure of dimension $k$ is 
only possible in the case that the underlying structure of dimension 
$k-1$ arises from a nilmanifold.

\end{document}